\documentclass[11pt]{article}
\usepackage{amssymb,amsmath,mathrsfs,amsfonts,graphicx,epsf} 
\usepackage{color}
\usepackage[dvipsnames]{xcolor}

 \input xy
 \xyoption{all}

 \usepackage{esint}

 \usepackage{fullpage}

\usepackage{authblk}

  \newcommand{\ffoot}[1]{}

\newcommand{\KKK}{{\cal K}}
\newcommand{\cut}{{\bf Cut}}

\newcommand{\bD}{\Delta}

\newcommand{\f}{\mathfrak{f}}

\usepackage{mathrsfs}
\usepackage{graphicx}
\usepackage{amssymb}
\usepackage{color}
\newtheorem{theorem}{Theorem}
\newtheorem{corollary}{Corollary}
\newtheorem{lemma}{Lemma}
\newtheorem{proposition}{Proposition}
\newtheorem{definition}{Definition}
\newtheorem{remark}{Remark}

\newcommand{\bt}{\begin{theorem}}
\newcommand{\et}{\end{theorem}}
\newcommand{\bl}{\begin{lemma}}
\newcommand{\el}{\end{lemma}}
\newcommand{\bp}{\begin{proposition}}
\newcommand{\ep}{\end{proposition}}
\newcommand{\bc}{\begin{corollary}}
\newcommand{\ec}{\end{corollary}}
\newcommand{\bdeff}{\begin{definition}}
\newcommand{\edeff}{\end{definition}}
\newcommand{\brem}{\begin{remark}}
\newcommand{\erem}{\end{remark}}
\renewcommand{\r}[1]{(\ref{#1})}
\newcommand{\con}{{\mathcal C}}

\newcommand{\bi}{\begin{itemize}}
\newcommand{\iii}{\item}
\newcommand{\ei}{\end{itemize}}
\newcommand{\bd}{\begin{description}}
\newcommand{\ed}{\end{description}}

\newcommand{\bqn}{\begin{eqnarray}}
\newcommand{\eqn}{\end{eqnarray}}
\newcommand{\eqnn}{\nonumber\end{eqnarray}}

\newcommand{\nn}{\nonumber}
\newcommand{\ba}[1]{\begin{array}{#1}}
\newcommand{\ea}{\end{array}}

\newcommand{\R}{\mathbb{R}}

\newcommand{\lam}{\lambda}
\newcommand{\g}{\gamma}
\newcommand{\al}{\alpha}
\newcommand{\eps}{\varepsilon}

\newcommand{\om}{\omega}

\newcommand{\VecM}{\mathrm{Vec}(M)}

\newcommand{\Gq}{{\gg}_q}

\newcommand{\Zz}{\mathcal{Z}}
\renewcommand{\gg}{{\bf G}}
\newcommand{\sign}{\mathrm{sign}}

\xdefinecolor{og}{named}{OliveGreen}

\newcommand{\e}{\mbox{e}}

\newcommand{\HH}{{\bf (H0)}}
\newcommand{\HHH}{{\bf (H1)}}

\newcommand{\frd}[2]{\frac{d #1}{d #2}}
\newcommand{\frp}[2]{\frac{\partial #1}{\partial #2}}

\def\F{{\bf{F}}}
\def\cn{\mathrm{cn\,}}
\def\sn{\mathrm{sn\,}}
\def\dn{\mathrm{dn\,}}

\def\sign{\mathrm{sign}}

\newcommand{\Oo}{O}
\newcommand{\Ss}{Q}

\newcommand{\tp}{${\mathbf T}^{\oplus}$}
\newcommand{\tm}{${\mathbf T}^{\ominus}$}
\newcommand{\ccc}{C\!\!\!\!C^\infty}
\newcommand{\set}{{\clubsuit}}

\title{\LARGE \bf 
Normal forms and invariants  for 2-dimensional almost-Riemannian structures\thanks{This research has been supported  by the European Research Council, ERC
StG 2009 ``GeCoMethods", contract number 239748, by the ANR ``GCM", program ``Blanc--CSD"
project number NT09-504490, and by the DIGITEO project ``CONGEO".}
}

\author[$\spadesuit$]{U.~Boscain} 
\affil[$\spadesuit$]{CNRS, CMAP \'Ecole Polytechnique, Palaiseau, France and
Team GECO, INRIA Saclay --
\^Ile-de-France {\tt ugo.boscain@cmap.polytechnique.fr}}
\author[$\clubsuit$]{G.~Charlot}
\affil[$\clubsuit$]{Institut Fourier, UMR 5582, CNRS/Universit\'e Grenoble 1, 100 rue des Maths, BP 74, 38402 St Martin d'H\`eres, France
and Team GECO, INRIA Saclay --
\^Ile-de-France 

{\tt Gregoire.Charlot@ujf-grenoble.fr}}
\author[$\dag$]{R.~Ghezzi} 
\affil[$\dag$]{CMAP \'Ecole Polytechnique, Palaiseau, France, 
Team GECO, INRIA Saclay --
\^Ile-de-France, 
and Department of Mathematical Sciences Rutgers University Camden NJ, 

{\tt roberta.ghezzi@rutgers.edu},
{\tt ghezzi@cmap.polytechnique.fr}}

\begin{document}

\maketitle

\begin{abstract}
\noindent
Two-dimensional almost-Riemannian structures are generalized Riemannian structures on surfaces for which a local orthonormal frame is given by a Lie bracket generating pair of vector fields that can become collinear. Generically, there are three types of points: Riemannian points where the two vector fields are linearly independent, Grushin points where the  two vector fields are collinear but their Lie bracket is not, and tangency points where the two vector fields and their Lie bracket are collinear and the missing direction is obtained with one more bracket.

In this paper we consider the problem of finding normal forms and functional invariants at each type of point. We also require that functional invariants are ``complete'' in the sense that they permit to recognize locally isometric structures.
 
The problem happens to be equivalent to the one of finding a smooth canonical parameterized curve passing through the point and being transversal to the distribution.

For  Riemannian points such that the gradient of the Gaussian curvature $K$ is different from zero, we use  the level set of $K$
 as support of the parameterized curve. 
For  Riemannian points such that the gradient of the curvature vanishes (and under additional generic conditions),  
we use a curve 
which is found by looking for crests and  valleys of the curvature. For Grushin points we use the set where the vector fields are parallel.

 Tangency points are the most complicated to deal with.
  The cut locus from the tangency point is not a good candidate as canonical parameterized curve since it is known to be non-smooth. Thus, we analyse the cut locus from the singular set and we prove that it is not smooth either.  A good candidate appears 
 to be a curve which is found by looking for crests and  valleys of the Gaussian curvature.  We prove that the support of such a curve is uniquely determined and has a canonical parametrization. 
\end{abstract}

\section{Introduction}\label{intro}

A $2$-dimensional Almost Riemannian  Structure ($2$-ARS for short) is a rank-varying sub-Riemannian structure  that can be  locally defined by a pair of smooth vector fields on a $2$-dimensional manifold, satisfying the H\"ormander condition  (see for instance \cite{AgrBarBoscbook,bellaiche,jean1,jean2}).   These vector fields play the role of an orthonormal frame. 
 It can also be defined as an Euclidean bundle of rank 2 on a 2-D manifold $M$ and a morphism of vector bundles from $E$ to $TM$ which gives rise to a Lie bracket generating distribution.

Let us denote by $\bD(q)$  the linear span of the two vector fields at a  point $q$.  Where  $\bD(q)$ is $2$-dimensional, the corresponding metric is Riemannian. Where $\bD(q)$ is $1$-dimensional, the corresponding Riemannian metric is not well-defined. However, thanks to the H\"ormander condition, one can still define the Carnot-Caratheodory distance between two points, which happens to be finite and continuous. 

$2$-ARSs were introduced in the context of hypoelliptic operators \cite{baouendi,FL1,grushin1}.  They appeared in  problems of  population transfer in quantum systems  \cite{q4,BCha,q1} and have applications to  orbital transfer in space mechanics  \cite{BCa,tannaka}.

Generically, the singular set  $\Zz$, where $\bD(q)$ has dimension $1$,  is a $1$-dimensional embedded submanifold  (see \cite{ABS}). There are three types of points: Riemannian points, Grushin points where   $\Delta(q)$  is $1$-dimensional and 
$\dim(\bD(q)+[\bD,\bD](q))=2$ and tangency points where $\dim(\bD(q)+[\bD,\bD](q))=1$ and the missing direction is obtained with one more bracket. One can easily show that at  Grushin points,
$\bD(q)$ is
 transversal to $\Zz$. Generically,  at tangency points  $\Delta(q)$ is tangent to $\Zz$ and  tangency points are isolated.

$2$-ARSs present very interesting phenomena. For instance, the presence of a singular set permits the conjugate locus  to be nonempty even if the Gaussian curvature is negative, where it is defined (see \cite{ABS}). Moreover, a Gauss--Bonnet-type formula can be obtained. More precisely, in 
 \cite{ABS,high-order} the authors studied the generic case without tangency points. In \cite{euler} this formula was generalized to the case in which tangency points are present. (For generalizations of Gauss--Bonnet formula in related contexts, see also \cite{agra-gauss, pelletier,pelle2}.) 
 In \cite{BCGS} a necessary and sufficient condition for two  2-ARSs  on the same compact manifold $M$ to be   Lipschitz equivalent was given. This equivalence was established  in terms of graphs associated with the structures. In \cite{camillo} the heat and the Schr\"odinger equation with the Laplace--Beltrami operator on a 2-ARS were studied. In that paper it was proven that the singular set acts as a barrier for the heat flow and for a quantum particle, even though geodesics can pass through the singular set without singularities.

 In this paper we consider the problem of  finding, at each type of point, a normal form which is completely reduced, in the sense that it depends only on the $2$-ARS and not on its local representation.\footnote{In this paper, an object that depends on the 2-ARS and not
on its local representation is called {\em canonical}.}
This consists in finding a canonical choice
for a local system of coordinates and for a local orthonormal frame, i.e., two vector fields $F_1$ and $F_2$ defined in a neighborhood of the origin on $\R^2$.\footnote{ To be able to fix completely the system of coordinates and the orthonormal frame,  and avoid the problem of having quantities defined up to sign, in this paper we assume that the 2-ARS is {\em totally oriented}, i.e., both the base manifold $M$ and the Euclidean bundle are oriented. See Definition \ref{d-oriented}.} Notice that, the classical normal form in the Riemannian case given by the exponential coordinates is
not completely reduced in our sense, since it is defined up to a choice of an orthonormal frame at the starting point. Rather, for the
completely reduced normal forms, which we are looking for in this paper, the coordinates and the orthonormal frame should be uniquely determined.

Once 
a canonical choice of $(F_1,F_2)$ is provided, the two components of  $F_1$ and the two components of $F_2$ are functional invariants of the structure, in the sense that locally isometric structures have the same components. Moreover, they are a {\em complete set of invariants} because they permits to recognize locally isometric structures: if two structures have the same invariants in a neighborhood of a point, then they are locally isometric.

However, one expects that some of these invariants are trivial in the sense that they have always the same value.
Let us make a rough computation of how many functional invariants we do expect. Among the 4 components of the two vector fields, we can fix two by choice of the coordinate system (which is defined by a diffeomorphism, i.e., by two functions of two variables) and one by choice of the orthonormal frame (which is defined by the choice of an angle, i.e., by a function of two variables). Indeed, one of the results of our paper is that under generic conditions, there is  a canonical way of choosing  an orthonormal frame $(F_1,F_2)$ and a system of coordinates  wherein  $F_1=(1,0)$ and $F_2=(0,f)$. Here $f$ is a complete set of invariants (that we call a {\em complete invariant} since it is just one function).

\begin{figure}[h!]
\begin{center}
\input{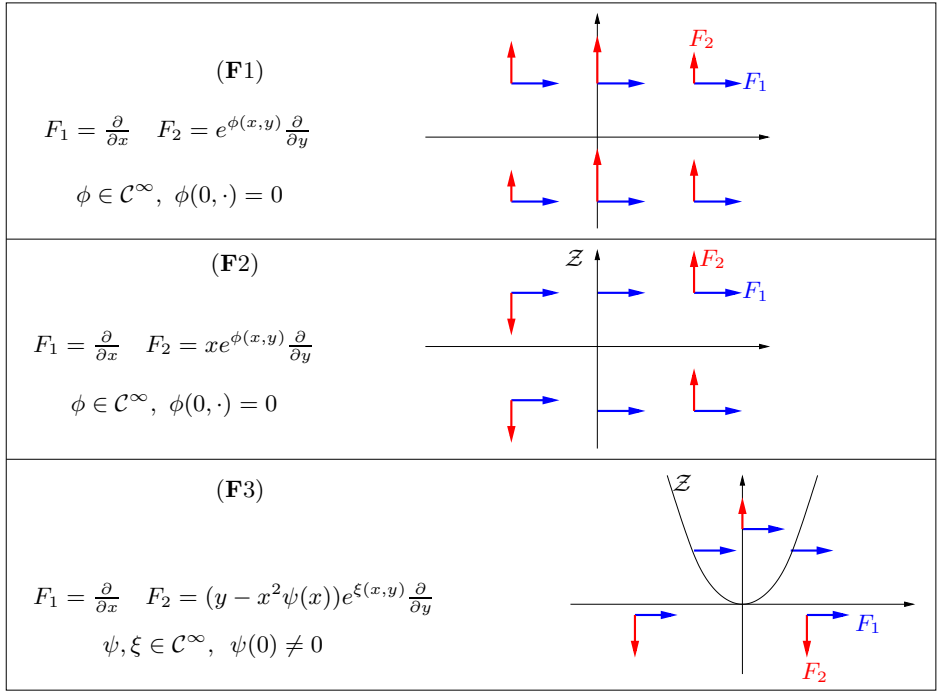}
\caption{The local representations established in \cite{ABS}}\label{fig-prenormal}
\end{center}
\end{figure}

Notice that the problem of finding a complete set of invariants is not completely trivial even in the simplest case of Riemannian points. See for instance the discussion in \cite{zelclass, kulkarni}.  Indeed even if one is able to canonically fix a system of coordinates, the Gaussian curvature in that system of coordinates is an invariant, but it is not a complete  invariant: there are non-locally isometric structures having the same curvature (an example is given in Section \ref{example}).  

A first step in finding normal forms has been realized in  \cite{ABS}, where  the local representations  given in Figure~\ref{fig-prenormal}  were found. However, the ones corresponding to Riemannian and tangency points  are not completely reduced. 
Indeed, there  exist changes of coordinates and rotations  of the frame for which an orthonormal basis has the same expression as in (\F1)  (resp. (\F3)), but with a different function $\phi$ (resp. with different functions $\psi$ and $\xi$).

\begin{figure}[h!]
\begin{center}
\input{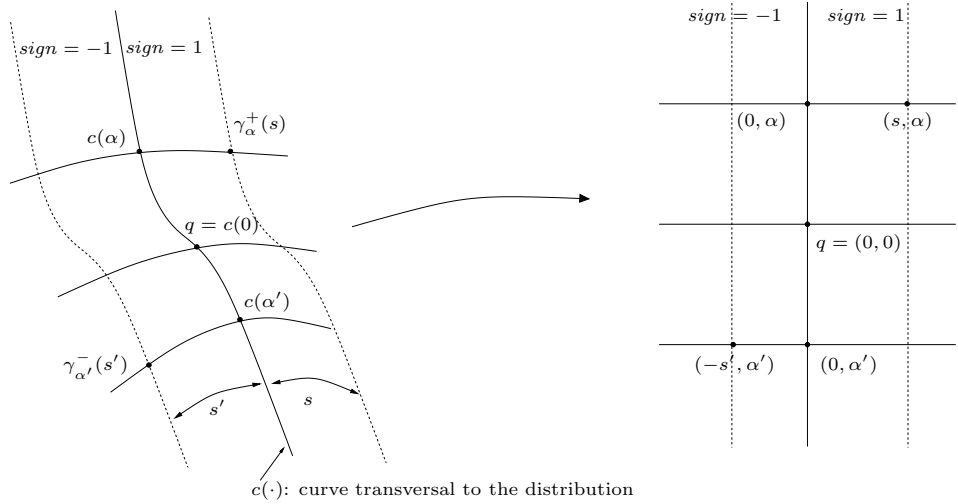}
\caption{The construction of coordinates starting from a parameterized curve $c(\cdot):]-\varepsilon,\varepsilon[\to M$. We denote by $\g^{\pm}_\alpha$ the geodesic starting at $c(\alpha)$, parameterized by arclength, entering the region where {\it sign} $=\pm 1$ and such that $d(\g^{\pm}_\alpha(s),c(]-\varepsilon,\varepsilon[))=s$.  As the distribution is transversal to $c(\cdot)$ at each point, the distance from $c(]-\varepsilon,\varepsilon[)$ is smooth.}\label{fig-costr}
\end{center}
\end{figure}

In order  to build  the coordinate system to which the  local expressions found in  \cite{ABS} apply,
the following idea  was  used. Consider a smooth parameterized curve passing through a point $q$.
If the curve is assumed to be transversal to the distribution 
at each point, then the Carnot--Caratheodory distance from the curve
 is shown to be smooth on a neighborhood of $q$  (see \cite{ABS}). Given
a point $p$ near $q$, the first coordinate of $p$ is, by definition, 
 the distance between $p$ and the chosen curve,  with a suitable choice of sign. The second coordinate of $p$ is the 
 parameter corresponding to the point (on the chosen curve) that realizes the distance between $p$ and the curve (see Figure~\ref{fig-costr}).
 If the parameterized curve used in this construction can be canonically built, then one gets 
 a local representation of the form $F_1=(1,0)$, $F_2=(0,f)$ which  cannot be further reduced. Hence, the functions $f$ is  a ``complete invariant'' in the sense above.
  
For Riemannian points, a canonical parameterized curve transversal to the distribution can be easily identified, at least at  points where the gradient of the Gaussian curvature is non-zero: one can use the level set of the curvature passing from the point, parameterized by arclength (see Section~ \ref{ss-transversal} and \ref{tcR1}). For points where the gradient of the curvature vanishes, under additional generic conditions,   we prove the existence of a smooth parameterized canonical curve passing through the point (a crest or a valley of the curvature, see Sections \ref{s-creste}, \ref{ss-transversal} and \ref{tcR2}). 
 
  For Grushin points, a canonical curve transversal to the distribution is the set $\Zz$. This curve has also a natural parameterization, as explained in Section \ref{ss-transversal}, and was used to get the local representation (\F2)  in Figure~\ref{fig-prenormal} (that, as a consequence, cannot be further reduced).
  
Concerning the local expression (\F3) in Figure~\ref{fig-prenormal},  in \cite{ABS} 
the choice of the smooth parameterized curve was arbitrary and not canonical. 
The main purpose of this paper is to find a canonical one.  Once this is done, one automatically gets a normal form which cannot be further reduced and the corresponding functional invariant at a tangency point.

\begin{figure}[h!]
\begin{center}
\includegraphics[width=6cm]{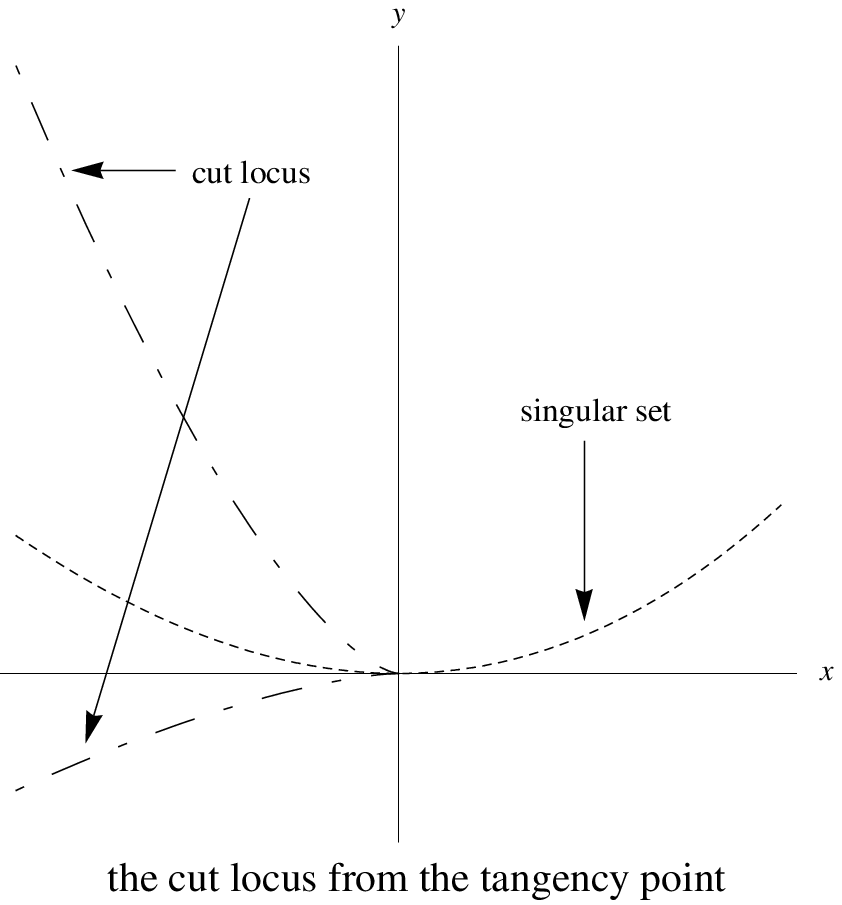}
\hspace{2cm}
\includegraphics[width=6cm]{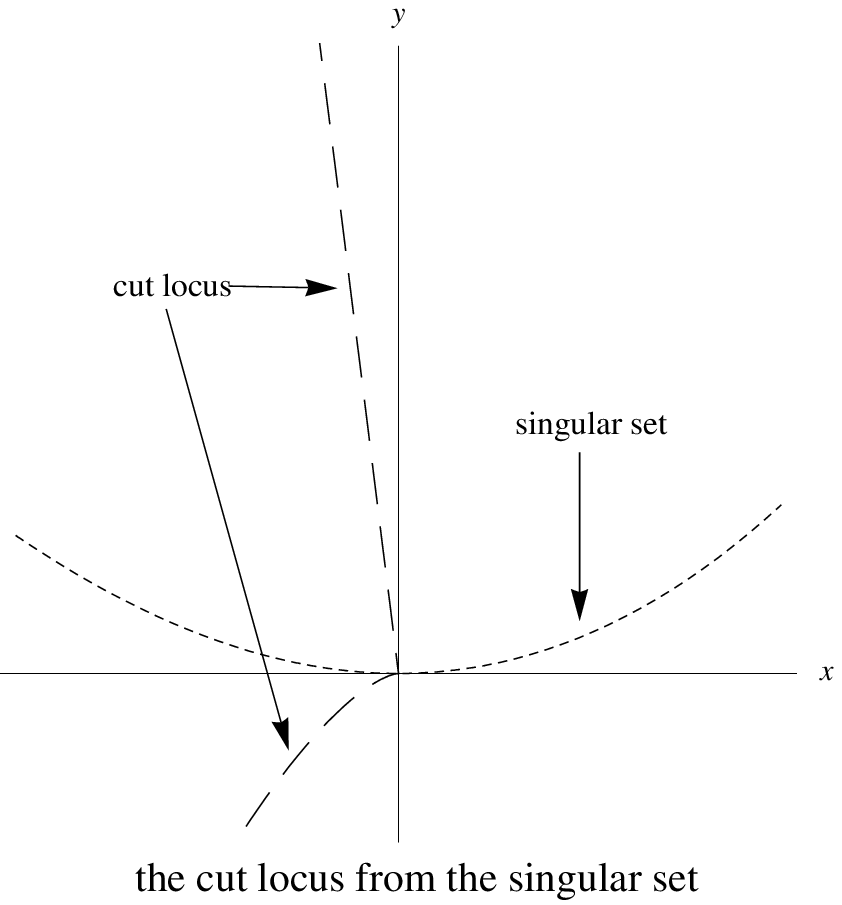}

\vspace{1cm}

\includegraphics[width=6cm]{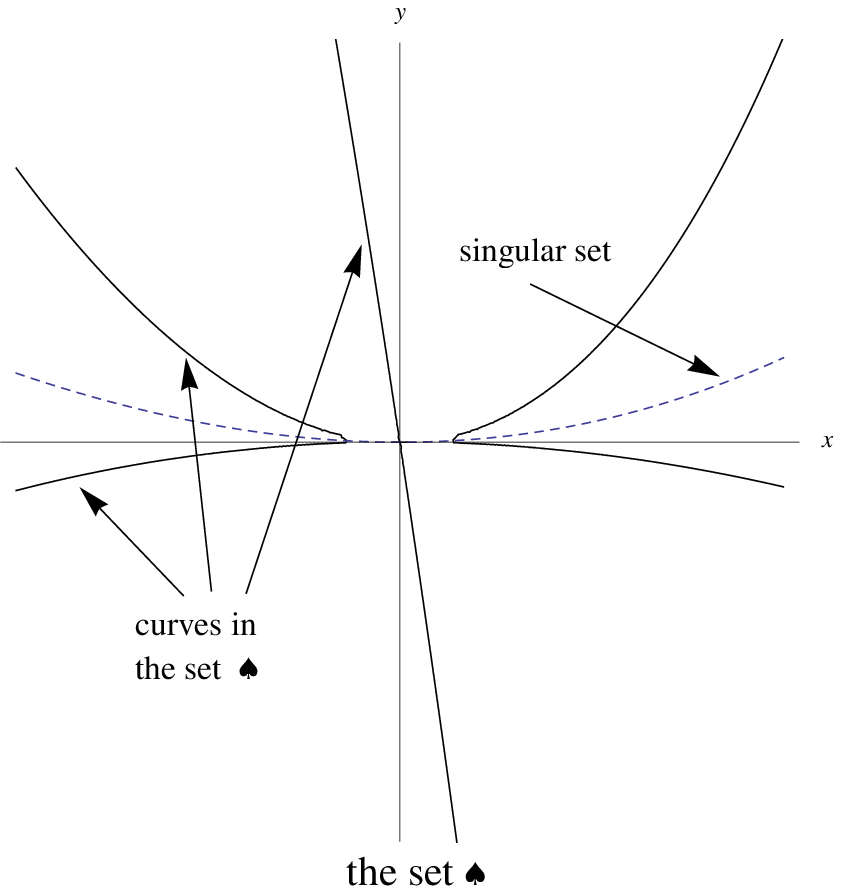}
\hspace{2cm}
\includegraphics[width=6cm]{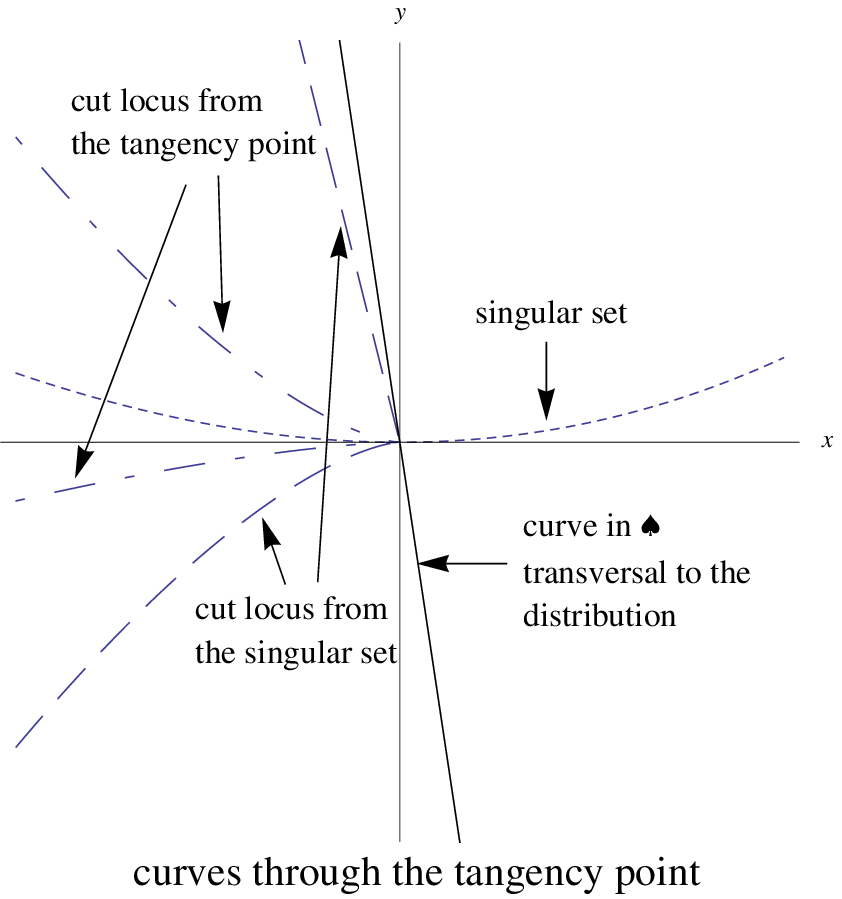}
\end{center}
\caption{  The singular set (dotted line), the cut locus from a tangency point 
(semidashed line),
the cut locus from the singular set (dashed line), and the set $\spadesuit$ (solid lines)  
for the ARS with orthonormal frame $F_1=\frp{}{x},\, F_2=(y-x^2-x^3)\frp{}
{y}$. In this case  the set $\spadesuit$ is the union of three smooth curves. Notice that all these curves but one are 
tangent to the distribution.}
\label{figurone}
\end{figure}

The most natural candidate for such a curve is the cut locus from the tangency point. Nevertheless,
this is not a good choice, as  in \cite{BCGJ} it was proven that in general   
the cut locus  from the point is not smooth but has an asymmetric cusp (see Figure  \ref{figurone}). 
Another possible candidate is the cut locus from the singular set in a neighborhood 
of the tangency point. The first result of the paper concerns the analysis  and the description of this locus: in Theorem~\ref{prop-cutz} we prove that
the cut locus from $\Zz$ is non-smooth in a neighborhood of a tangency point  (see Figure \ref{figurone} for an example). Even if not useful for the construction of the completely reduced normal form, this result is a step forward in understanding the geometry of tangency points.

\begin{figure}[h!]
\begin{center}
\includegraphics[width=12truecm]{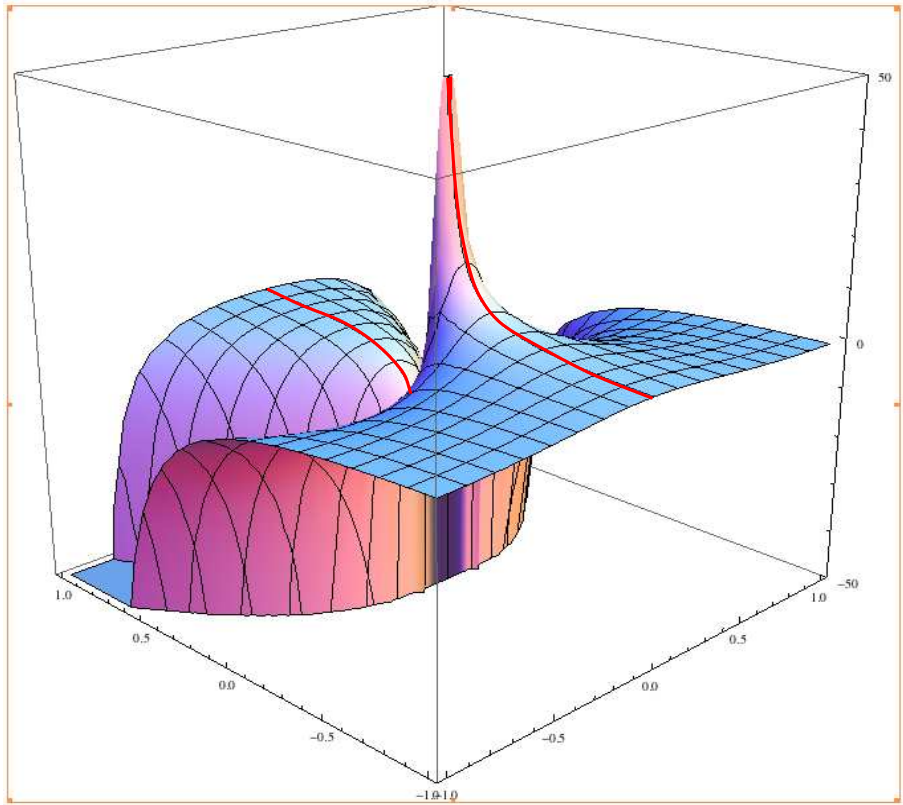}
\caption{A crest of the curvature at a tangency point for the almost-Riemannian structure on $\R^2$ having 
$X(x,y)=(1,0),\, Y(x,y)=(0, y-x^2)$ as orthonormal frame.}\label{fig-KRESTACCIA}
\end{center}
\end{figure}

A third possibility is to look for curves which are crests or
valleys of the Gaussian curvature and  intersect transversally the singular set at a tangency point.
The second result of the paper (see Theorem~\ref{rnf-tangency})  consists in the  proof of  the existence of  such a curve.  Moreover, this curve admits a canonical regular  parameterization.
 An example of a crest of the curvature at a tangency point is shown in Figure~\ref{fig-KRESTACCIA}.

Notice that  for 2-ARSs, tangency points are the most difficult to handle 
due to the fact that the asymptotic of the distance to the singular set is different from the two sides of the singular set.  
In \cite{BCGJ} the authors gave a description of the geometry of the nilpotent approximation at a tangency point. Also, they provided jets of the exponential map and a description of the cut and conjugate loci from a tangency point in the generic case. 

However, tangency points are far to be deeply understood.  An open question is the convergence or the divergence of the integral of the geodesic curvature on the boundary of  a tubular  neighborhood  of the singular set, close to a tangency point. This question arose in the proof of the Gauss--Bonnet theorem given in \cite{euler}. In that paper, thanks to numerical simulations, the authors conjecture  the divergence of such integral.

The structure of the paper is the following.  In Section~\ref{sec-basdef}, we briefly recall the notion of almost-Riemannian structure, we define the concept of local representations and describe the procedure to build a local representation  from a parameterized curve transversal to the distribution. Finally, we  define the set containing crests and valleys of the curvature. 

Section~\ref{s-CRNF} is devoted to the formal definition of completely reduced normal forms and of complete set of invariants as maps that to a germ of a 2-ARS associate a choice of local orthonormal frame and a choice of a set of functions which permits to distinguish local isometric structures.

In Section~\ref{s-main}, we state the main results: the construction for each type of point  of a canonical parameterized curve transversal to the distribution, the description of the cut locus from the singular set at a tangency point and the construction of the normal forms and of the invariants. 

The next sections are devoted to the proof of the results.  Section~\ref{sec-cutz} is devoted to the description of  the cut locus from the singular set in a neighborhood of a tangency point. In Section  \ref{s-6}, we prove the existence of a canonical parameterized curve passing transversally to the distribution at a tangency point. In Section \ref{s-7}, we prove Corollary \ref{main-theorem} which describes the completely reduced normal form and the complete invariant.

 \section{Preliminaries}\label{sec-basdef}

In this section we recall some basic definitions  in the framework
 of  2-ARS following  \cite{ABS,euler}.

Let $M$ be a  smooth  connected surface without boundary.  Throughout the paper, unless specified,  manifolds are smooth (i.e., $\con^{\infty}$) and without boundary. Vector fields  and differential forms are smooth. The set of smooth vector fields
 on $M$ is denoted by $\VecM$.

 \begin{definition}
 A \emph{ $2$-dimensional almost-Riemannian structure} (2-ARS) is a triple
${\mathcal S}=
(E, \f,\langle\cdot,\cdot\rangle)$
 where 
  $E$ is a  vector bundle of rank $2$ over $M$ and $\langle\cdot,\cdot\rangle$ is an Euclidean structure on $E$, that is, $\langle\cdot,\cdot\rangle_q$ is a   scalar product 
on $E_q$ smoothly depending on $q$.  Finally
   $\f:E\rightarrow TM$ is a morphism of vector bundles,    
i.e., {\bf (i)}  the diagram 
$$
\xymatrix{
 E  \ar[r]^{\f} \ar[dr]_{\pi_E}   & TM \ar[d]^{\pi}            \\
 & M                          
}    
$$
 commutes, where  $\pi:TM\rightarrow M$ and $\pi_E:E\rightarrow M$ denote
  the canonical projections and {\bf (ii)} $\f$ is linear on fibers.
Denoting by $\Gamma(E)$ the $\con^\infty(M)$-module of 
smooth sections on $E$,  we assume the submodule $\bD=\{\f\circ\sigma\mid\sigma\in\Gamma(E)\}$ to be bracket generating, i.e.,
  $Lie_q(\bD)= T_qM$ for every $q\in M$.
\end{definition}

Denote by $\Zz$ the singular set of $\f$, i.e., the set of points $q$ of $M$ such that $\mbox{dim}(\f(E_q))$ is less than 2.

\bdeff
\label{d-oriented}
A 2-ARS is said to be {\em oriented} if $E$ is oriented as vector bundle. We say that a 2-ARS is {\em totally oriented} if both $E$ and $M$ are oriented. For a totally oriented 2-ARS, $M$ is split into two open sets $M^+$, $M^-$ such that $\Zz=\partial M^+=\partial M^-$, $\f:E|_{M^+}\rightarrow TM^+$ is an orientation-preserving isomorphism and $\f:E|_{M^-}\rightarrow TM^-$ is an orientation reversing-isomorphism.
\edeff

A property $(P)$ defined for 2-ARSs  is said to be {\it generic}
if for every rank-2 vector bundle $E$   over $M$, $(P)$ holds for every $\f$ in an  open and dense  subset of the set of  morphisms of vector bundles from $E$ to $TM$, endowed with the 
$\con^\infty$-Whitney topology.

Let ${\mathcal S}=(E,\f,\langle\cdot,\cdot\rangle)$ be a  2-ARS on a surface $M$.
We  denote by $\bD(q)$  the linear subspace $\{V(q)\mid  V\in \bD\}=\f(E_q)\subseteq  T_q M$. 
 The Euclidean structure on $E$ induces a symmetric 
positive-definite bilinear form $G:\bD\times\bD\to\con^\infty(M)$ defined by
 $G(V,W)=\langle\sigma_V,\sigma_W\rangle$ where $\sigma_V,\sigma_W$ are the unique sections of $E$ satisfying $\f\circ\sigma_V=V,\, \f\circ\sigma_W=W$. At points $q\in M$ where $\f|_{E_q}$ is an isomorphism, $G$ is a tensor and the value $G(V,W)|_q$ depends only on $(V(q), W(q))$.  This is no longer true at points  $q$ where $\f|_{E_q}$ is not injective.

If $(\sigma_1,\sigma_2)$ is an orthonormal frame for $\langle\cdot,\cdot\rangle$ on
 an open subset $\Omega$ of $M$, an {\it  orthonormal frame for  ${\mathcal S}$} on $\Omega$ is 
the pair  $(\f\circ\sigma_1,\f\circ\sigma_2)$. 

For every $q\in M$  
and every $v\in\bD(q)$ define
$
\Gq(v)=\inf\{\langle u, u\rangle_q \mid u\in E_q,\f(u)=v\}.
$

An  absolutely continuous curve $\g:[0,T]\to M$ 
 is  {\it admissible} for ${\mathcal S}$ 
if  
there exists a measurable essentially bounded function 
$[0,T]\ni t\mapsto u(t)\in E_{\g(t)}
$ such that 
$\dot \g(t)=\f(u(t))$  for almost every $t\in[0,T]$. 
Given an admissible 
curve $\g:[0,T]\to M$, the {\it length of $\g$} is  
\bqn
\ell(\g)= \int_{0}^{T} \sqrt{ \gg_{\gamma(t)}(\dot \g(t))}~dt.
\eqnn
The {\it Carnot-Caratheodory distance} (or sub-Riemannian distance) on $M$  associated with 
${\mathcal S}$ is defined as
\bqn\nonumber
d(q_0,q_1)=\inf \{\ell(\g)\mid \g(0)=q_0,\g(T)=q_1, \g\ \mathrm{admissible}\}.
\eqn
The finiteness and the continuity of $d(\cdot,\cdot)$ with respect 
to the topology of $M$ are guaranteed by  the Lie bracket generating 
assumption  (see \cite{book2}).  
The Carnot-Caratheodory distance    endows $M$ with the 
structure of metric space compatible with the topology of $M$ as differentiable manifold.

Locally, the problem of finding a curve realizing the distance between two fixed points  $q_0,q_1\in M$ is naturally formulated as the distributional optimal control problem 
\bqn
\dot q=\sum_{i=1}^2 u_i F_i(q)\,,~~~u_i\in\R\,,
~~~\int_0^T 
\sqrt{\sum_{i=1}^2 u_i^2(t)}~dt\to\min,~~q(0)=q_0,~~~q(T)=q_1,\eqnn
where $F_1,F_2$ is a local orthonormal frame for the structure.

A {\it geodesic} for  ${\cal S}$  is an admissible 
curve $\g:[0,T]\to M$, such that 
 $\gg_{\gamma(t)}(\dot \g(t))$ is constant and   
for every sufficiently small interval 
$[t_1,t_2]\subset [0,T]$, $\g|_{[t_1,t_2]}$ is a minimizer of $\ell$. 
A geodesic for which $\gg_{\gamma(t)}(\dot \g(t))$ is (constantly) 
equal to one is said to be parameterized by arclength. 

 If  $(F_1,F_2)$ is an orthonormal frame on an open set $\Omega$, a curve parameterized by arclength is a geodesic if and only if it is the projection on $\Omega$ of a solution of the Hamiltonian system corresponding to the Hamiltonian
\bqn\label{eq:hamiltonian}
H(q,{\bf p})=\frac12( ({\bf p}F_1(q))^2 + ({\bf p}F_2(q))^2),~~q\in\Omega,\,{\bf p}\in T^\ast_q\Omega.
\eqn
lying on the level set  $H=1/2$. This follows from the Pontryagin Maximum Principle \cite{pontryagin-book} in the case of 2-ARS. Its simple form follows from the absence of abnormal extremals in 2-ARS, as a consequence of the H\"ormander condition, see \cite{ABS}. Notice that $H$ is well-defined on the entire $T^*M$,  since  formula \r{eq:hamiltonian} does not depend on the choice of the orthonormal frame.
 When looking for a geodesic $\gamma$  realizing the distance from  a submanifold  $N$ (possibly of dimension zero), one should add the transversality condition ${\bf p}(0)T_{\gamma(0)}N=0$.

The cut locus $\cut_N$ from  $N$ is the set of points $p$ 
for which there exists a geodesic realizing the distance between $N$ and $p$ losing optimality after $p$.
 It is well known (see for instance \cite{agrompatto} for a proof in the three-dimensional contact case)
that, when there are no abnormal extremals, if $p\in \cut_{N}$ then one of the following two possibilities occour: {\bf i)} more than one minimizing geodesic  reaches $p$; {\bf ii)} $p$ belongs to the first conjugate locus from $N$ defined as follows. To simplify the notation, assume that all geodesics are defined on $[0,\infty[$.
Define
\bqn
C^0=\{\lambda=(q,{\bf p})\in T^*M\mid q\in N,~H(q,{\bf p})=1/2,~{\bf p} T_{q}N=0\}
\eqnn
and 
\bqn
&&\exp:C^0\times[0,\infty[\to M\nn\\
&&~~~~~~~~(\lam,t)\mapsto\pi(e^{t \vec{H}}\lam)
\eqnn
where $\pi$ is the canonical projection $(q,{\bf p})\to q$ and $\vec{H}$ is the Hamiltonian vector field corresponding to $H$. The first conjugate time from $N$ for the geodesic $\exp(\lambda,\cdot)$ is
\bqn
t_{conj}(\lam) = \min\{t > 0, (\lam,t) \mbox{ is a critical point of }\exp\}.
\eqnn
and the first conjugate locus from $N$ is $\{\exp(\lam,t_{conj}(\lam))\mid  \lam\in C^0 \}$.

\subsection{Local representations}\label{sec-nf}

     Set $\bD_1 = \bD$ and $\bD_{k+1}=\bD_k+[\bD,\bD_k]$.   Let us introduce the main
     assumptions under which all the results of the paper are proven. 
 
\bd
\iii[\HH] {(i)} $\Zz$ is an
embedded one-dimensional 
submanifold of
$M$;

\hspace{.1cm}{(ii)} the points $q\in M$ where $\bD_2(q)$ is
one-dimensional are isolated;

\hspace{.1cm}{(iii)}  $\bD_3(q)=T_qM$ for every $q\in M$.
\ed
Property \HH\ is generic for 2-ARSs (see \cite{ABS}).

      \bdeff\label{def:locrep}
A  \emph{local representation} of a 2-ARS at  a point $q\in M$ is a pair of vector fields  $(X,Y)$ on $\R^2$ such that there exist: {\bf i) }a neighborhood $U$ of $q$ in $M$, a neighborhood $V$ of $(0,0)$ in $\R^2$ and a diffeomorphism $\varphi:U\rightarrow V$ such that $\varphi(q)=(0,0)$; {\bf ii)} a local orthonormal frame $(F_1,F_2)$ of $\Delta$ around $q$, such that  $\varphi_*F_1=X$, $\varphi_*F_2=Y$, where $\varphi_*$ denotes the push-forward.
\edeff

%

Let us state a result which will be crucial in the following.

\bp[\cite{ABS}]\label{localreprcan}
Under the hypothesis \HH, it is always possible to get a local representation under the form $X=\partial_x$, $Y=f(x,y)\partial_y$, where $f$ is a smooth function such that one of the following conditions holds:
$f(0,0)\neq0$,
$\partial_xf(0,0)\neq0$,
$\partial_{xx}f(0,0)\neq0$.
\ep

In the following, we give a procedure which permits to build a local representation of the form $(\partial_x, f(x,y)\partial_y)$ starting from
a totally oriented 2-ARS and a parameterized curve transversal to the distribution. 
This procedure provides a completely reduced normal form once a canonical transversal curve is identified and 
a proof of Proposition \ref{localreprcan}.\footnote{Proposition \ref{localreprcan} is formulated for
non totally oriented 2-ARS. One can use Procedure 1 to prove it by fixing arbitrary local orientations on the manifold and on the Euclidean bundle.}

\medskip

\noindent{\bf Procedure 1:}
\begin{enumerate}
\iii Choose any smooth parametrized curve $c(\cdot):]-\varepsilon,\varepsilon[\to M$ such that $c(0)=q$, $\dot{c}(\alpha)\neq 0$ for
$\alpha\in]-\varepsilon,\varepsilon[$ and 
$\mbox{span}(\dot c(\alpha))+\bD(c(\alpha))=T_{c(\alpha)}M$ for every $\alpha\in]-\eps,\eps[$.

\iii Denote by ${\bf p}:]-\varepsilon,\varepsilon[\to T^*M$ the smooth map such that, for all $\alpha\in]-\varepsilon,\varepsilon[$,
\begin{enumerate}
\item ${\bf p}(\alpha)\in T^*_{c(\alpha)}M$,
\item ${\bf p}(\alpha)\dot c(\alpha)=0$,
\item $H(c(\alpha),{\bf p}(\alpha))=\frac 1 2$, where $H$ is defined in formula (\ref{eq:hamiltonian}),
\item ${\bf p}(\alpha)V(\al)>0$ (for every $\al\in]-\eps,\eps[$), if $(V(\al),\dot c(\alpha))$ is a positively oriented (with respect to the orientation of $M$) pair of vectors applied in $c(\al)$.
\end{enumerate}
 Remark that the map ${\bf p}$ is unique once $c(\cdot)$ is fixed.
\iii Define ${\cal E}:\R^2\to M$ as the map that associates with the pair $(x,y)$ the projection on $M$ 
of the solution at time $x$ of the hamiltonian system on $T^*M$ associated with $H$ with initial condition 
$(c(y),{\bf p}(y))$. The map ${\cal E}$ is a local diffeomorphism from a neighborhood $V$ of $0\in\R^2$ to
a neighborhood $U$ of $q\in M$ and preserves the orientation. Define $\varphi$ as the inverse of the restriction of ${\cal E}$ to these neighborhoods.
\iii Let $\sigma\in\Gamma(E|_U)$ be the unique section of norm one such that $\varphi_*(\f\circ\sigma)=\partial_x$ 
and $\rho$ such that $(\sigma,\rho)\in\Gamma(E|_U)^2$ is a positively oriented orthonormal frame of $E|_U$. 
\iii Define $X$ and $Y$ in $\mathrm{Vec }(\R^2)$ by $X=\varphi_*(\f\circ\sigma)$ 
and $Y=\varphi_*(\f\circ\rho)$.
\iii Define $f:V\subset\R^2\to \R$ by $Y(x,y)=f(x,y)\partial_y$.
\end{enumerate}

\bl\label{completed} Procedure 1 can be completed.\el

\noindent {\bf Proof:} Step 1 is possible thanks to the fact that $\Delta(q)$ has dimension at least one.

Step 2 is possible by simple linear algebra considerations. This choice of ${\bf p}$ corresponds to the initial condition
of a geodesic (solution of the Hamiltonian system defined by $H$) of the 2-ARS, transversal to the 
curve $c(\cdot)$ at the point $c(\alpha)$, hence minimizing locally the distance to the support of $c(\cdot)$.

In Step 3, the fact that ${\cal E}$ is a local diffeomorphism holds true since
$\frac{\partial {\cal E}}{\partial x}(0)$, $\frac{\partial {\cal E}}{\partial y}(0)$ are linearly independent vectors,
which implies that the Jacobian of ${\cal E}$ at $0$ is not $0$. Indeed $\frac{\partial {\cal E}}{\partial x}(0)$
is the initial velocity of a geodesic transversal to the curve $c(\cdot)$ when $\frac{\partial {\cal E}}{\partial y}(0)=\dot c (0)$.
As a consequence $\varphi$ is a local coordinate system in the neighborhood of $q$.

In Step 4 the existence of $\sigma$ is guaranteed by the fact that $x\mapsto{\cal E}(x,y)$ are geodesics for all $y\in]-\varepsilon, \varepsilon[$
which implies that ${\cal E}_*\partial_x$ is in the distribution $\Delta$. Once $\sigma$ is defined, the existence of $\rho$ is a trivial fact.

In Step 5 the existence of $X$ and $Y$ is guaranteed by the fact that $\varphi$ is a local diffeomorphism.

In Step 6, the fact that $Y$ can be written as $f(x,y)\partial_y$ is the consequence of the Pontryagin Maximum Principle. 
Indeed, since the curves $x\mapsto {\cal E}(x,y)$ are geodesics minimizing the distance to the support of $c(\cdot)$ and since these curves
have the form $x\mapsto (x,y)$ in the $(x,y)$-coordinates, then the curves $x\mapsto(x,y)$ realize the distance between vertical lines.
This implies that any solution of the Pontryagin Maximum Principle with initial condition transversal to the vertical axis  $\{(x,y)\mid x=0\}$ should be transversal 
to the vertical lines for each $x$. Hence it  should annihilate $\partial_y$. As a consequence of the Pontryagin Maximum Principle it should also annihilate
the orthogonal to $\partial_x$ for the metric. Hence the orthogonal to the metric is proportional to $\partial_y$. The fact that
$f$ is smooth is just a consequence of the smoothness of the metric and of ${\cal E}$.\hfill$\blacksquare$

\brem
Notice that Procedure 1 provides
\begin{itemize}
\iii a local coordinate system $\varphi$ around $q$ such that $\varphi(q)=(0,0)$ and $\varphi$ preserves the orientation of $M$;
\iii a 2-ARS on a neighborhood of $0\in\R^2$ having as  positively oriented orthonormal frame 
 $(\partial_x,f(x,y)\partial_y)$.
 One can easily check that, thanks to \HH, $f$ satisfies at least one of the following conditions $f(0,0)\neq0$,
$\partial_xf(0,0)\neq0$, 
$\partial_{xx}f(0,0)\neq0$.
\end{itemize}
\erem

A consequence of Lemma \ref{completed} is the following.
\begin{corollary} In the totally oriented case, constructing a 
local representation of the form $X=\partial_x$, $Y=f(x,y)\partial_y$ is equivalent to choose 
a parameterized curve transversal to the distribution.
\end{corollary}

In \cite{ABS} the following possible forms for the function $f$ in Proposition~\ref{localreprcan}  were found. 
\bp[\cite{ABS}]
If a 2-ARS satisfies \HH, then
 for every point $q\in M$ there exist a local representation having one of the forms
$$
\begin{array}{lll}
(\F1) & F_1(x,y)=\frp{}{x}, & F_2(x,y)=e^{\phi(x,y)}\frp{}{y},\\
(\F2) & F_1(x,y)=\frp{}{x}, & F_2(x,y)=xe^{\phi(x,y)}\frp{}{y},\\
(\F3) & F_1(x,y)=\frp{}{x}, & F_2(x,y)=(y-x^2\psi(x))e^{\xi(x,y)}\frp{}{y},
\end{array}
$$
where $\phi$, $\psi$ and $\xi$ are smooth functions such that $\phi(0,y)=0$ and $\psi(0)\neq0$.
\label{prop-prenormal}\ep

\brem
Notice that in the proposition above we do not take  account of orientations.
\erem
\bdeff Under hypothesis \HH,
a point $q$ is said to be a {\em Riemannian} point   if $\bD(q) = T_qM$,   {\em Grushin} point
if $\bD(q)$ is one-dimensional and $\bD_2(q)=T_qM$, {\em tangency} point if $\bD(q)=\bD_2(q)$ is $1$-dimensional and $\bD_3(q)=T_qM$.
 \edeff
Local representations for a Riemannian, Grushin and  tangency points are given by  (\F1), (\F2), (\F3) respectively.

 \brem
 The local representations {\em (\F1)} and  {\em (\F3)} were obtained   arbitrarily choosing a parameterized curve transversal to the distribution. Hence, they are not completely reduced  in the sense that there are isometric structures 
 having  local representations of type {\em (\F1)}  (resp.~{\em(\F3)}), but with  different functions $\phi$ (resp.~with different functions $\psi$ and $\xi$). In other words, the functions appearing in   {\em(\F1)}  and   {\em(\F3)} are not invariants of the structure. 
 
This topic, together with the proof that the local representation of type  {\em (\F2)} is completely reduced,  is discussed  starting from Section \ref{s-CRNF}.
 \erem

 \subsection{Crests and valleys of the curvature}
\label{s-creste}
\newcommand{\kk}{K}

 In this section we find an equation whose solutions contain ``crests'' and ``valleys'' of the curvature. These loci are  used in Section \ref{s-main} to construct a canonical parameterized curve transversal to the distribution and to obtain invariants of the structure.

For simplicity let us assume that the curvature $\kk$ is a Morse function on $M\setminus\Zz$, so that its level sets are locally either one-dimensional manifolds (at regular points) or isolated points (at maxima or minima) or the union of two one-dimensional manifolds which transversally intersect  (at saddle points). 

\begin{figure}[h!]
\begin{center}
\input{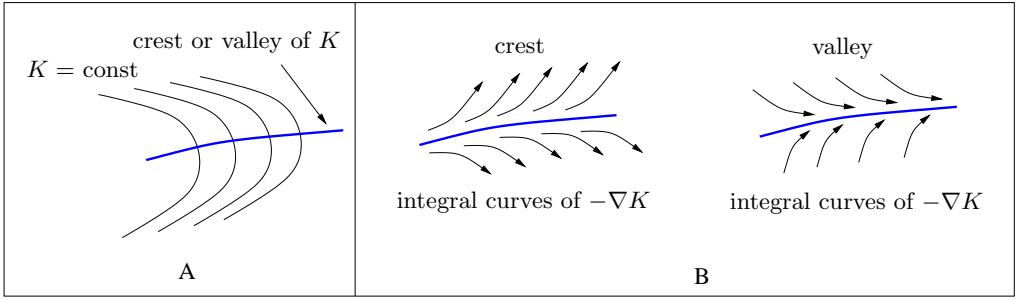}
\caption{Crests and valleys of the curvature}\label{fig-creste-valli}
\end{center}
\end{figure}

Let us think to $\kk$ as a function describing the altitude of a mountain. Roughly speaking,  crests and  valleys of $\kk$ are the loci where the distance among the level sets of $\kk$ has a local maximum along a level set of $K$, see Figure \ref{fig-creste-valli}A.
To distinguish among crests and valleys one should consider integral curves of  $-\nabla \kk$.   Here $\nabla \kk$ denotes
 the almost-Riemannian gradient of $\kk$, i.e.,  the unique vector  such that $G(\nabla K,\cdot)=dK(\cdot)$.
Locally,  integral curves of $-\nabla \kk$ diverge from a crest and converge to a valley, see Figure \ref{fig-creste-valli}B.

To find an equation satisfied by crests and valleys, let $C$ be a level set of $\kk$. For simplicity, assume  that $C$ is a one-dimensional manifold at $q$. The point $q$ belongs to a crest or a valley of $\kk$ if $\|\nabla \kk\|^2
:=G(\nabla \kk,\nabla \kk)$ has a local minimum on $C$ at $q$. 
A necessary condition for this to happen is that $\nabla\|\nabla \kk\|^2$ and $\nabla \kk$ are collinear at $q$. Namely, crests and valleys are loci lying in the set
\bqn
\label{pre-picche}
\{q\in M\mid G(\nabla\|\nabla \kk\|^2,\nabla \kk^\perp )|_q=0\},
\eqn
where $(\nabla K)^\perp$ is a non-vanishing vector satisfying
  $G(\nabla K, (\nabla K)^\perp)=0$.
Notice that this set contains not only crests and valleys but also {\it anti-crests} (i.e., loci in which $\|\nabla \kk\|^2$ has a local maximum along a level set of $\kk$ and integral curves of $-\nabla \kk$ diverge from these loci)   and {\it anti-valleys} (i.e., loci in which  $\|\nabla \kk\|^2$ has a local maximum along a level set of $\kk$ and integral curves of $-\nabla \kk$ converge to these loci). 
Moreover the set \r{pre-picche} contains critical points of $\kk$ and $\|\nabla \kk\|^2$.

For 2-ARSs, $K$ is not  defined on the singular set $\Zz$, in particular it is not defined at tangency points. Since below we look for a canonical curve passing through a tangency point, it is  convenient to look for it inside the set $\spadesuit$ defined in the following. 

\bdeff\label{def:picca}
Let ${\mathcal S}=
(E, \f,\langle\cdot,\cdot\rangle)$ be a 2-ARS. Define
$$\spadesuit:= \{q\in M\setminus \Zz\mid  G(\nabla\|\nabla \kk\|^2,\nabla \kk^\perp )|_q=0\}\cup\{q\in M\mid q \mbox{ is a tangency point of } {\cal S}\}.$$ 
\edeff

\brem
Notice that, for a 2-ARS defined by an orthonormal frame of the form $(\partial_x,f(x,y)\partial_y)$, at Riemannian points the metric and the curvature are given by
(see for instance \cite{ABS})
\bqn\label{curvature}
g=
\left( \ba{cc} 
1&0\\0&\frac{1}{f^2}
\ea\right),~~~
K=\frac{-2(\partial_x f)^2+f\partial^2_xf}{f^2}.
\eqn
\erem

\section{Completely reduced normal forms and invariants}

\label{s-CRNF}

In the following, $2$-ARSs are assumed to be totally oriented unless specified.
 Recall that  under the generic conditions \HH\ we have three types of points: Riemannian ({\bf R} for short), Grushin  ({\bf G} for short) and tangency points  ({\bf T} for short). 

In the following, we need to distinguish two types of tangency points.

\begin{figure}[h!]
\begin{center}
\input{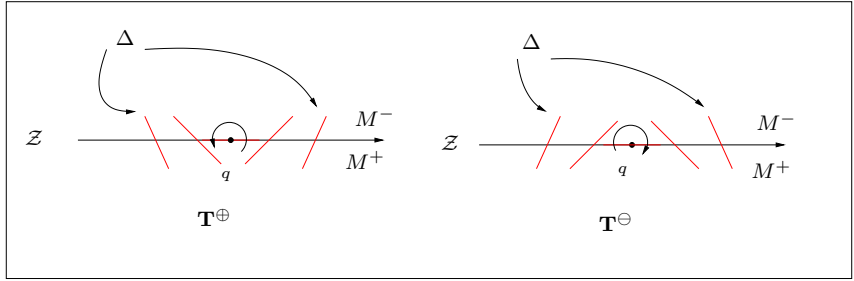}
\caption{ The two different types of tangency points\label{tauu}}
\end{center}
\end{figure}

\bdeff\label{deftp}
Let $q$ be a tangency point. Let us orient $\Zz$ as boundary of $M^-$ (see Figure \ref{tauu}).
 We say that the tangency point is of type \tp\  (resp. \tm) if the distribution is rotating positively (resp. negatively) along $\Zz$ at $q$. See Figure \ref{tauu}.
\edeff

\brem 
Notice that the distinction between \tp\ and \tm\ was used also in \cite{euler,BCGS},  to obtain respectively a Gauss-Bonnet Theorem for  2-ARSs and a classification of 2-ARSs w.r.t.~Lipschitz equivalence (\tp\  corresponds to a contribution $\tau_q=-1$ and \tm\  corresponds to a contribution $\tau_q=1$, where $\tau_q$ is defined in \cite{euler}).
\erem 

In the following we have to make  additional generic conditions to treat Riemannian points. 

\bd
\iii{HA.} On $M\setminus\Zz$ the Gaussian curvature is a Morse function (i.e., at points where 
       $\nabla K=0$,  the Hessian is non-degenerate). 
        Moreover, if $q$ is such that  $\nabla K=0$, then the Hessian of $K$ at $q$ computed in an orthonormal system of coordinates has two distinct eigenvalues.   
\ed   
 It is a standard fact that this condition is generic. It ensures that $\nabla K$ vanishes only at isolated points and that at these points 
the curvature has  {\bf i)} a local minimum, {\bf ii)} a local maximum or {\bf iii)} a saddle.  In case {\bf iii)}  the curvature has one crest and one valley intersecting transversally.\footnote{Here for convenience we consider that if $\nabla K(q)=0$ then $q$ belongs to the crests, valleys, anticrests, antivalleys reaching the point.}
The condition on the eigenvalues of the Hessian implies that: in case
 {\bf i)} the curvature has one valley and one anti-valley intersecting transversally;  in case {\bf ii)}  the curvature has one crest and one anti-crest intersecting transversally.

To be able to fix an orientation on crests and valleys we need a higher order condition.
\bd
\iii{HB.}  Assume HA. If $c(\cdot):]-\eps,\eps[\to M$ is a smooth curve parameterized by arclength the support of which is contained in a crest or a valley of $K$ and such that $\nabla K(c(0))=0$ then $\frac{\partial^3 }{\partial t^3}K(c(t))|_{t=0}\neq0$.
\ed   

In the following we call \HHH\ the collection of  HA and HB. We need to treat separately two types of Riemannian points. 

\bdeff
Assume \HHH\ and let $q\in M\setminus\Zz$. We say that  $q$ is a Riemannian point of type {\bf R1} if $\nabla K(q)\neq0$ and of  type {\bf R2} if $\nabla K(q)=0$.
\edeff

Recall that a pair of vector fields $(X_1,X_2)$ defined in a neighborhood  of the origin of $\R^2$ for which dim$(Lie_0\{X_1,X_2\})=2$ defines a totally oriented 2-ARS, namely the one for which $(X_1,X_2)$  is a positively oriented orthonormal frame and the orientation on $\R^2$ is the canonical one.

 \bdeff   
Let  $\Ss$ be the set of germs of  totally oriented 2-ARS verifying \HH\ and \HHH. 

Let $\Oo$ be the set of germs at the origin of pair of vector fields on $\R^2$ 
such that if $ o\in \Oo$ then 
  {\bf (i)} $o$ is Lie bracket generating, {\bf (ii)} the corresponding germ ${\cal S}_o$ of totally oriented 2-ARS belongs to $\Ss$. 
   
We say that two germs of totally oriented 2-ARS are {\em p-isometric} if they are isometric and the isometry preserves 
the orientation of the base manifolds and of the vector bundles.\footnote{A local isometry between the base manifolds preserves the orientation on the vector bundles if the pushforward of a positively oriented orthonormal frame for the first structure is a positively oriented orthonormal frame for the second structure.}
  \edeff

 \bdeff[Completely reduced normal form] A {\em completely reduced normal form} (CRNF, for short) for  totally oriented 2-ARSs is a map $N:\Ss\to\Oo$  which associates  with a germ ${\cal S}$ of a 2-ARS the germ $o=N({\cal S})$ at the origin of a pair of vector fields on $\R^2$ such that
 \bi
 \iii  ${\cal S}$ and ${\cal S}_o$ are  p-isometric;
 \iii  $N({\cal S}_1)=N({\cal S}_2)$ if and only if ${\cal S}_1$ and ${\cal S}_2$ are p-isometric.
 \ei

  \edeff

 \bdeff[Complete set of invariants]
Let ${\ccc}$ be the set of germs at the origin of $\con^\infty$ functions on $\R^2$. 
 A {\em complete set  of invariants of cardinality} $k$ for totally oriented 2-ARSs  is a map  $I:\Ss\to \left(\ccc\right)^k$ such that $I({\cal S}_1)=I({\cal S}_2)$ if
 and only if ${\cal S}_1$ and ${\cal S}_2$ are p-isometric. When the complete set of invariants is of cardinality one, it is called a {\em complete invariant}.
  \edeff
 
  \brem
Notice that once a CRNF is obtained, one can build a complete set of invariants for totally oriented 2-ARSs by constructing it on $\Oo$, since each element of $\Ss$ is p-isometric to an element of $\Oo$. One can simply take
 the non trivial functions among the 4 components of the two vector fields of $o\in\Oo$
 as a complete set of invariants. 
  \erem

\subsection{On 2-d Riemannian manifolds the curvature is not a complete invariant } \label{example}

On 2-d Riemannian manifolds, even when a canonical system of coordinates can be constructed (for instance, via a canonical transversal curve and following Procedure 1), the curvature written in this system of coordinates is not a complete invariant.
Consider for instance on $]-\frac 1 2,\frac 1 2[\times\R$ the two Riemannian metrics,
 $$g_1=\left(\ba{cc}
1 & 0\\ 0 & \frac 1 {(x+1)^2}
\ea\right),\;\;\;\; g_2=\left(\ba{cc}
1 & 0\\ 0 & (x+1)^4
\ea\right).$$

Both these metrics are written in a system of coordinates built with Procedure 1, by taking as parametrized curve the level set of $K$
passing through the origin and parametrized by arclength. From (\ref{curvature}), one gets for both metrics that
$$
K=-\frac 2 {(x+1)^2}.
$$
However, one easily proves that these two metrics are not isometric, even locally.

 \subsection{A CRNF by using a canonical parameterized curve transversal to the distribution}
 
To build a CRNF,  we find, for each type of point {\bf R1}, {\bf R2}, {\bf G}, {\bf T}, a canonical parameterized curve transversal to the distribution and  we use  Procedure 1 to build a local representation.
These curves are built in Section \ref{ss-transversal}.

\medskip
 Let $\Gamma$ be the set of germs of smooth curves taking value on a 2-dimensional manifold. For a point $q$ of a 2-d manifold, let us denote $\Gamma_q$ the subset of germs $c(\cdot)\in\Gamma$ such that $c(0)=q$.   
We have to build a map $\chi:\Ss\ni{\cal S}\mapsto c(\cdot)\in\Gamma$ such that: (i)  $c(\cdot)\in\Gamma_q$, 
where $q$ is the base point of ${\cal S}$; (ii) $c(\cdot)$ is transversal to the distribution at $q$; (iii) $\chi$ is invariant by p-isometries, that is, if $\phi$
is a p-isometry between ${\cal S}$ and ${\cal S}'$ then $\phi(\chi({\cal S}))=\chi({\cal S}')$.

Once $\chi$ is built, by using  Procedure 1, one can obtain a CRNF $N_\chi:\Ss\to\Oo$.
By Proposition \ref{localreprcan} we have  the following fact.
 \bp
The image of $N_\chi$ is a subset of  $\Oo$ of elements of the form $(X_1,X_2)$ such that $X_1=(1,0)$ and $X_2=(0,f)$,
where $f$ is a germ of a smooth function on $\R^2$.
 \ep

 Since $\chi$ is invariant with respect to p-isometries, $N_\chi$ provides automatically a complete invariant.
 
\bp
 Let  $I_\chi:\Ss\to \con^\infty(\R^2,\R)$ be the map $I_\chi({\cal S})=f$ where $f$ is the function appearing as second component of  the second vector field of  $N_\chi({\cal S})$. Then $I_\chi$  is a complete invariant.
 \ep
 
The next section is devoted to the construction of the map $\chi$. Namely, for each type of point  we build a parameterized curve transversal to the distribution which is invariant by p-isometries. We build the map $I_\chi$ and study its image.  
 
 \section{Main results}
 \label{s-main}
 \subsection{Looking for a transversal curve}
 \label{ss-transversal}
 In this section we build explicitly the map $\chi$ which associates with the germ of a totally oriented 2-ARS at the point $q$  a parameterized curve transversal to the distribution at $q$. By construction this curve is canonical. We treat  the different types of points {\bf G}, {\bf T}, {\bf R1}, {\bf R2}, separately.  

\begin{figure}
\begin{center}
\input{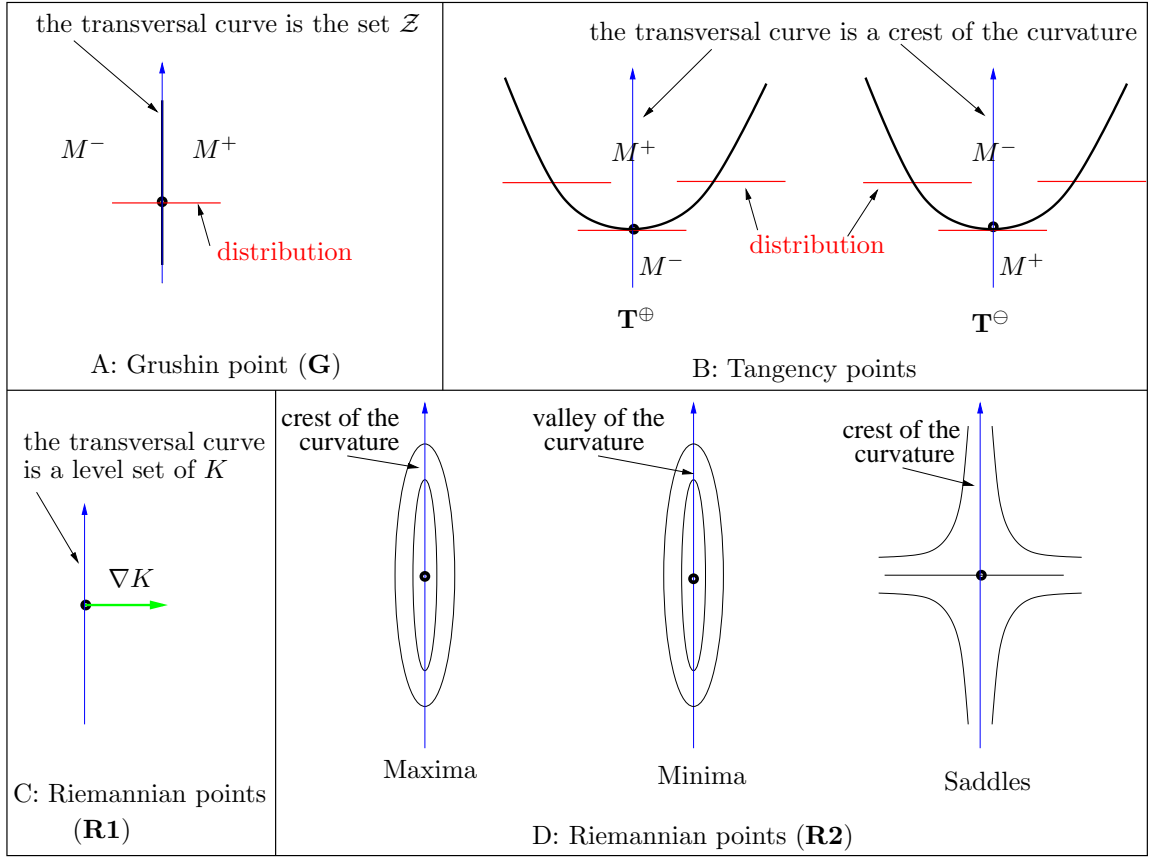}
\caption{Canonical parameterized transversal curves at the different type of points.\label{f-TUTTIPUNTI}}
\end{center}
\end{figure}
%

\subsubsection{Grushin points}
 Let $q$ be a   Grushin point.
In a neighborhood of $q$,  $\Zz$ is transversal to the distribution. It is easy to see that   for any positively oriented orthonormal frame $(G_1,G_2)$, the Lie bracket $[G_1,G_2]|_{\Zz}$ modulo elements in $\bD$ does not change, and it is not zero close to $q$.
Hence we fix the transversal curve to be 
the parameterized curve $c(\cdot)$ having  the singular set $\Zz$
as support, such that $c(0)=q$ and $[G_1,G_2]|_{c(\alpha)}=  c'(\alpha)\mbox{ mod }\Delta$.
Notice that with this choice,  the orientation of  $c(\cdot)$ is such that 
on its right (w.r.t. the orientation of $M$) lies the set $M^+$, see Figure \ref{f-TUTTIPUNTI}A.
This construction is the same as the one made in  \cite{ABS} to build (\F2),  except that this time we take  account of orientations.

Notice that another possible candidate for  a transversal curve could be the cut locus from $q$. However,  this choice would require to prove that the cut locus from $q$ is smooth.

\subsubsection{Tangency points}
  
  The case of tangency points is rather complicated. The first candidate as support of a smooth curve is the cut locus from the tangency point. Let us recall a result of  \cite{BCGJ} where the shape of the cut locus 
at a tangency point has been computed.
\bp[\cite{BCGJ}]\label{formacut}
Let $q\in M$ be a tangency point of a 2-ARS satisfying  \HH\ and assume there exists a local representation of the type {\em (\F3)} at $q$ with the property
$$
\psi'(0)+\psi(0)\partial_x\xi(0,0)\neq0.
$$
Then the cut locus from the tangency point accumulates at   $q$ as an
asymmetric cusp whose branches are  locally separated by $\Zz$.
In the coordinate system where the chosen local representation is
{\em (\F3)}, the cut locus is locally
$$
\{(\sign(\al_1)t^2,\sqrt{|\al_1|} t^3+o(t^3))\mid t>0\}\cup
\{(\sign(\al_2)t^2,-\sqrt{|\al_2|} t^3+o(t^3))\mid t>0\},
$$
with
$\al_i=c_i/(\psi'(0)+\psi(0)\partial_x\xi(0,0))^3$, the constants
$c_i$ being nonzero and independent on the structure.
\ep
In general, due to Proposition~\ref{formacut}, the cut locus from $q$ is neither smooth nor transversal to the distribution at $q$.

Another candidate would be the cut locus from $\Zz$ in a neighborhood of a tangency point. A description of such locus is given by the following theorem.

\bt\label{prop-cutz}
Let $q\in M$ be a tangency point of a 2-ARS satisfying  \HH\ and assume there exists a local representation of the type {\em (\F3)} at $q$ with the property
$$
\alpha=\psi'(0)+\psi(0)\partial_x\xi(0,0)\neq0.
$$
Then the cut locus from the singular set $\Zz$ in a neighborhood of $q$ accumulates at   $q$ as the union of two curves locally separated  by $\Zz$. 
 One of them is contained in the set $\{y>x^2\psi(x)\}$,   takes the form
$$
\{(-1/2\psi'(0)t+o(t),t+o(t))\mid t>0\},
$$
and accumulates at $q$ transversally to the distribution.
The other one is contained in the set $\{y<x^2\psi(x)\}$ and takes the form
$$
\{(\alpha \, \om \, t^2+o(t^2),- t^3+o(t^3))\mid t>0\},
$$
where $\om\neq 0 $ is a constant depending on the structure. This part of the cut locus accumulates at $q$ 
with
 tangent direction at $q$ belonging to the distribution. 
%
%
\et
The proof of Theorem~\ref{prop-cutz}    
is given in Section \ref{sec-cutz}.

As a consequence of the Theorem~\ref{prop-cutz}, in general the cut locus from $\Zz$ is neither smooth nor transversal to the distribution at $q$.

Finally, we  look for a smooth curve, transversal to the distribution, the support of which is a subset of $\spadesuit$ (see Definition~\ref{def:picca}).
More precisely, we get the following result proved in Section \ref{s-6}.

\bt\label{rnf-tangency} Let ${\cal S}$ be a totally oriented 2-ARS on $M$ satisfying \HH\ and let $q\in M$ be a
tangency point. 
There exist  $\epsilon>0$ and a unique smooth parameterized curve $c(\cdot)$ defined on $]-\epsilon,\epsilon[$ which satisfies the following properties: 
(i) $c(0)=q$, $\dot c(0)\notin \bD(q)$; 
(ii) the support of $c(\cdot)$ is contained in $\spadesuit$; 
(iii) if  $(X,Y)$ is the positively oriented orthonormal frame constructed following Procedure 1, there exists $\lambda>0$ such that  
\bi
\iii $[X,[X,Y]]|_{(0,y)}=-2\partial_y$ and $[\partial_y,Y]|_{(0,0)}=\lambda \partial_y$ if $q$ is of type \tp. 
\iii $[X,[X,Y]]|_{(0,y)}=2\partial_y$ and $[\partial_y,Y]|_{(0,0)}=-\lambda \partial_y$ if $q$ is of type  \tm. 
\ei
\et
Notice that with this choice,    $\dot c(0)$ points towards $M^+$ for \tp\ and towards  $M^-$ for \tm, see Figure \ref{f-TUTTIPUNTI}B. It is always in the direction where the cut locus from $\Zz$ accumulates to the tangency point transversally to $\Zz$ (see Theorem~\ref{prop-cutz}).

\subsubsection{Riemannian points of type 1}\label{tcR1}
 Let $q$ be a point such that $\nabla K(q)\neq0$.  There exists a unique arclength parameterized curve $c(\cdot)$ such that $c(0)=q$, which is supported in the level set $\{p\in M~|~K(p)=K(q)\}$, and along which the gradient of $K$ points on the right of $c(\cdot)$  (with respect to the orientation of $M$). See Figure \ref{f-TUTTIPUNTI}C.
 
\subsubsection{Riemannian points of type 2}\label{tcR2}

Let $q$ be a point such that $\nabla K(q)=0$. Thanks to \HHH, at maxima and saddles, there exists exactly one smooth curve $\gamma_c(\cdot)$ 
parameterized by arclength with $\gamma_c(0)=q$, whose support is a crest, and such that 
$\frac{\partial^3 }{\partial t^3}K(\gamma_c(t))|_{t=0}>0$. Similarly, at minima and saddles, there exists exactly 
one smooth curve $\gamma_v(\cdot)$  parameterized by arclength  with $\gamma_v(0)=q$,  whose support is a valley, and such that 
$\frac{\partial^3 }{\partial t^3}K(\gamma_v(t))|_{t=0}>0$. 
When $q$ is a maximum or a saddle, we set $c(\cdot)=\gamma_c(\cdot)$. If $q$ is a minimum we set $c(\cdot)=\gamma_v(\cdot)$.
Notice that the support of $c(\cdot)$ is contained in $\spadesuit$. See Figure  \ref{f-TUTTIPUNTI}D.

\subsection{Completely reduced normal forms and invariants}

In this section, using the canonical parameterized curves built above and applying Procedure 1, we give a CRNF for generic totally oriented 2-ARS which provides a complete invariant, i.e., a map $I_\chi$ which associates with a germ of  a 2-ARS the germ of a smooth function. Finally we characterize the image of $I_{\chi}$. Recall that $\ccc$ is the set of germs of smooth functions on $\R^2$.

 \bc [Main Result] \label{main-theorem} Let $N_{\chi}$ be the completely reduced normal form for totally oriented 2-ARS obtained from the map $\chi$ with  Procedure  1. Then the image of $N_{\chi}$ has the form $(X_1,X_2)$ where $X_1=(1,0)$, $X_2=(0,f(x,y))$ with $f$ satisfying the following properties.
\bd
  \iii[{[{\bf G}]}]  If $q$ is a point of type {\bf G} then $f$ satisfies
\bd
\iii {\em(Ga)} $f(0,\cdot)=0$.
\iii{\em (Gb)} $\partial_xf(0,\cdot)=1$.
\ed
In this case we write  $f\in \set_{\bf G}\subset\ccc$.
  \iii[{[{\bf T}]}] If $q$   is a point of type \tp\ (resp. \tm) then $f$ satisfies
  \bd
  \iii{\em (Ta)} $f(0,0)=0$.
  \iii{\em (Tb)} $\partial_xf(0,0)=0$.
  \iii{\em (Tc)} $\partial^2_xf(0,\cdot)=-2$ for  \tp  (resp. $\partial^2_xf(0,\cdot)=2$ for  \tm).
  \iii{\em (Td)} $\partial_yf(0,0)>0$   for  \tp (resp. $\partial_yf(0,0)<0$ for  \tm).
  \iii{\em (Te)} for $y\neq 0$, we have $G(\nabla\|\nabla \kk\|^2,\nabla \kk^\perp )|_{(0,y)}=0$ which ensures that the vertical axis is included in $\spadesuit$.
\ed
In this case we write  $f\in \set_{\oplus}\subset\ccc$ (resp. $f\in \set_{\ominus}$).
 \iii[{[{\bf ${\bf R}_1$}]}]  If $q$ is a point of type {\em {\bf R1}}  then $f$ satisfies
\bd
  \iii{\em (R1a)} $f(0,\cdot)=1$ (resp $-1$) if at $q$ the manifold $M$ and the oriented 2-ARS have the same (resp. opposite) orientation. 
  \iii{\em (R1b)} The second compontent of $\nabla K$ vanishes along the vertical axes. 
  \iii{\em (R1c)} The first component of $\nabla K$ is positive  along the vertical axes. 
\ed
In this case we write $f\in\set_{\bf R_1}\subset\ccc$. Writing $f=\pm e^{\phi(x,y)}$ these conditions read
\bd
  \iii{\em (R1d)} $\phi(0,\cdot)=0$.
  \iii{\em (R1e)} $-2\partial_x^2\phi(0,y)\partial_x\partial_y\phi(0,y)
+\partial_x^2\partial_y\phi(0,y)=0,$ for all $y$.
  \iii{\em (R1f)} $\partial_x^3\phi (0,y)-2 \partial_x\phi (0,y) \partial_x^2\phi(0,y)>0$, for all $y$.
\ed
 
 \iii[{[{\bf ${\bf R}_2$}]}]   If $q$ is a point of type {\em {\bf R2}}  then $f$ satisfies
\bd
\iii{\em (R2a)} $f(0,\cdot)=1$ (resp $-1$) if at $q$ the manifold $M$ and the oriented 2-ARS have the same (resp. opposite) orientation. 
\iii{\em (R2b)} For $y\neq 0$, we have $G(\nabla\|\nabla \kk\|^2,\nabla \kk^\perp )|_{(0,y)}=0$ which ensures that the vertical axis is included in $\spadesuit$.
 \iii{\em (R2c)} If $q$ is a local maximum for $K$, then $0>\partial^2_yK(0,0)>\partial^2_xK(0,0)$
which ensures that the vertical axis is a crest (and the horizontal one an anticrest).
 \iii{\em (R2d)}  If $q$ is a local minimum for $K$, then $0<\partial^2_yK(0,0)<\partial^2_xK(0,0)$
which ensures that the vertical axis is a valley (and the horizontal one an antivalley).
 \iii{\em (R2e)} If $q$ is a saddle for $K$, then $\partial^2_yK(0,0)>0>\partial^2_xK(0,0)$
which ensures that the vertical axis is a crest (and the horizontal one a valley).
 \iii{\em (R2f)} $\partial^3_yK(0,0)>0$ which fixes the orientation of the vertical axis.
\ed
In this case we write  $f\in\set_{\bf R_2}\subset\ccc$.
\ed
\ec
 
\brem
  Notice that the sets $\set_{\bf G}$,  $\set_\oplus$,  $\set_\ominus$,  $\set_{\bf R_1}$, $\set_{\bf R_2}$ are disjoint.

 Moreover their union is not $\ccc$. This is a consequence of the fact that we have generic conditions and 
  that $f$ is constructed using a canonical transversal curve.
When the function $f$ is obtained applying Procedure 1 to any
transversal curve,  it does not satisfy the conditions given in Corollary \ref{main-theorem} in general.
  \erem

\section{Proof of Theorem~\ref{prop-cutz}: the cut locus $\cut_\Zz$ from the singular set $\Zz$}\label{sec-cutz}

In this section, we prove Theorem~\ref{prop-cutz} starting from the local representation  (\F3). 
Notice that by applying the coordinates change 
$$
\tilde x= x,\quad \tilde y =\frac{y}{\psi(0)},
$$
we may assume that $\psi(0)=1$. For sake of readability, in the following we rename $\tilde x, \tilde y$ by $x, y$. Since 
$\psi(0)>0$,  the singular set $\Zz$ is locally contained in the upper half plane $\{(x,y)\mid y\geq 0\}$. 

Locally, the singular set separates $M$ in two domains $\{(x,y)\mid y-x^2\psi(x)>0\}$ and $\{(x,y)\mid y-x^2\psi(x)<0\}$. 
First, notice that  $\cut_\Zz\cap\Zz=\emptyset$, since we are computing the cut locus from $\Zz$. 
Second, thanks to hypothesis $\HH$, the only points of $\Zz$ where $\cut_\Zz$ may accumulate are the tangency points, 
since all other points of $\Zz$ are Grushin points, where $\Delta$ is transversal to $\Zz$. Hence, close to a tangency point, 
$\cut_\Zz$ is the union of two parts,  $\cut_\Zz^+$  lying in the upper domain $\{(x,y)\mid y-x^2\psi(x)>0\}$ and  
$\cut_\Zz^-$  in the lower one.

Applying the Pontryagin Maximum Principle, geodesics for the ARS  are projections on $\R^2$ of solutions of the Hamiltonian system 
associated with the function
$$
H=\frac{1}{2}(p_x^2+p_y^2(y-x^2\psi(x))^2\e^{2\xi(x,y)}),
$$
that is, solutions of the  system
\bqn\label{ham-sys}
\left\{
\ba{lll}
\dot x & = & p_x\\
\dot y & = & p_y((y-x^2\psi(x))\e^{\xi(x,y)})^2\\
\dot p_x & = & p_y^2(y-x^2\psi(x))(2x\psi(x)+x^2\psi'(x)-(y-x^2\psi(x))\frp{\xi}{x}(x,y))\e^{2\xi(x,y)}\\
\dot p_y & = & -p_y^2(y-x^2\psi(x))(1+(y-x^2\psi(x))\frp{\xi}{y}(x,y))\e^{2\xi(x,y)}.\\
\ea\right.
\eqn
In addition, a geodesic starting from $\Zz$ with $x(0)= a\neq 0$, realizing the distance from $\Zz$ and parameterized
by arclength, must satisfy the transversality condition 
$$
p_x(0)=\pm 1,\;\;\;p_y(0)=\mp \frac{1}{2a\psi(a)+a^2\psi'(a)}.
$$

As we shall see, the two components of $\cut_\Zz$ have different natures. In the upper domain the geodesic 
starting at a point $(a,a^2\psi(a))$ and transversal to $\Zz$ reaches its cut point at a time of order 1 in $|a|$. 
Whereas in the lower domain the geodesic starting at the same point reaches its cut point at a time of order 1 in $\sqrt{|a|}$.

Introducing  the new time  variable $s=\frac{t}{\eta}$ where $\eta>0$ is a parameter, system \r{ham-sys}  becomes
\bqn\label{ham-sys-rep}
\left\{
\ba{lll}
\frd{x}{s} & = & \eta p_x\\
\frd{y}{s} & = & \eta p_y((y-x^2\psi(x))\e^{\xi(x,y)})^2\\
\frd{p_x}{s} & = & \eta p_y^2(y-x^2\psi(x))(2x\psi(x)+x^2\psi'(x)-(y-x^2\psi(x))\frp{\xi}{x}(x,y))\e^{2\xi(x,y)}\\
\frd{p_y}{s} & = & -\eta p_y^2(y-x^2\psi(x))(1+(y-x^2\psi(x))\frp{\xi}{y}(x,y))\e^{2\xi(x,y)}.\\
\ea\right.
\eqn
In order to make the coordinates $x$, $y$, $p_x$ and $p_y$ dependent on the time variable and on the initial condition
$a$, 
in the following  we write them as functions of $t$ and $a$.  In Section~\ref{subs-sopra}, the parameter $p_x(0)$ is assumed to be 
$-\sign(a)$,  since this implies that the geodesic enters the upper domain. In Section~\ref{subs-sotto} it is assumed to be 
$\sign(a)$  since this implies that the geodesic enters the lower domain.

In the following, since we are studying $\cut_\Zz$ close to the tangency point $(0,0)$, $a$ can be assumed  as small as we want.
\subsection{The upper part of the cut locus}\label{subs-sopra}

We consider the geodesic starting from a point of $\Zz$, realizing for small time the distance from $\Zz$ and entering the upper domain, 
that is, with the initial conditions 
\bqn\label{destra-su}
\ba{ll}\displaystyle
x(t=0,a)=a , &  p_x(t=0,a)=-\sign(a),\\
y(t=0,a)= a^2\psi(a), & p_y(t=0,a)=\dfrac{\sign(a)}{2a\psi(a)+a^2\psi'(a)}.
\ea
\eqn

\noindent{\bf Computation of jets.} For $a>0$,
choosing $\eta=a$ in \r{ham-sys-rep} and writing $x$, $y$, $p_x$ and $p_y$ as functions of $a$ and $s$, one can check 
that if $x,y,p_x,p_y$ have orders  $1,2,0,-1$ in $a$ respectively, then the dynamics has the same or higher orders. As a consequence, 
since the initial conditions respect these orders, we can compute jets with respect to $a$ of the solution of system \r{ham-sys-rep} 
under the form
$$
\begin{array}{lllccrrl}
x(s,a) & = & a x_0(s) + a^2 x_1(s) +a^3\bar x(s,a) & & & p_x(s,a), & = & {p_x}_0(s) + a {p_x}_1(s) +a^2\bar p_x(s,a),\\
y(s,a) & = & a^2 y_0(s) + a^3 y_1(s) +a^4\bar y(s,a)& & & p_y(s,a), & = & a^{-1} {p_y}_0(s) + {p_y}_1(s) +a \bar p_y(s,a),
\end{array}
$$
where $\bar x, \bar y, \bar p_x, \bar p_y$ are smooth functions.
Using \r{destra-su}, the initial conditions are given by
$$
\begin{array}{llll}
x_0(0)=	1, & x_1(0)=0, & {p_x}_0(0)=-1, & {p_x}_1(0)=0,\\
y_0(0)=1, &y_1(0)=\psi'(0), & {p_y}_0(0)=\frac{1}{2}, & {p_y}_1(0)=-\frac{3}{4}\psi'(0),
\end{array}
$$
and from system \r{ham-sys-rep} we easily get 
$$
x_0(s)=1-s, \;\; x_1(s)\equiv0,\;\; y_0(s)\equiv1,\;\; y_1(s)\equiv\psi'(0),
$$
whence
$$
x(t,a)=a-t+a^3\bar x(a,t/a),\;\; y(t,a)=a^2+a^3\psi'(0)+a^4\bar y(a,t/a).
$$
Similarly, for $a <0$ one gets
$$
x(t,a)=a+t+a^3\bar x(a,t/a),\;\; y(t,a)= a^2+ a^3\psi'(0)+ a^4\bar y(a,t/a).
$$

\bl\label{derivee}
For $a$ small enough and $t$ such that $|\frac t a|<2$ one gets that
$$
\frp{x(t,a)}{a}>\frac 1 2, \quad \quad \frp{y(t,a)}{a}<0 \mbox{ if } a<0, \quad \quad \frp{y(t,a)}{a}>0 \mbox{ if } a>0.
$$
\el
The proof is a direct consequence of the computation of the jets of the geodesics.

\bl\label{calcolo-cut}
A geodesic with the parameter $a$ small enough intersects the geodesic with initial condition
$\bar a =-a-a^2\psi'(0)+o(a^2)$ at time $t_{int}(a)=|a|(1 +\frac 1 2 a\psi'(0)+o(a))$. The pair
$(\bar a, t_{int}(a))$ is unique among the pairs $(b,\tau)$ realizing the intersection and satisfying
$ab<0$, $\tau>0$, $0<|\frac{\tau}{a}|<2$ and $0<|\frac{\tau}{b}|<2$.
\el
A direct consequence of this lemma is that the cut time of a geodesic starting from $\Zz$ with initial condition $a$ is 
bounded from above by $t_{int}(a)$.

\noindent {\bf Proof:}  Assume $a>0$, the proof being the same for $a<0$.
Let $\bar a<0$ and $t>0$ be such that $0<\frac{t}{a}<2$ and $0<|\frac{t}{\bar{a}}|<2$. 
If $y(t,a)=y(t,\bar a)$ then
$$
a^2+a^3\psi'(0)+o(a^3)={\bar a}^2+{\bar a}^3\psi'(0)+o({\bar a}^3),
$$
whence $\bar a =-a -a^2\psi'(0)+o(a^2)$. Moreover, $x(t,a)=x(t,\bar a)$ implies
$t=a+\frac{1}{2}a^2\psi'(0)+o(a^2)$ and the intersection 
point is
\bqn\label{cut-sopra}
x_{int}(a)=-\frac{\psi'(0)}{2} a^2+o(a^2),\;\; y_{int}(a)=a^2+o(a^2).
\eqn
This, together with the fact that $\frac{\partial x(t,a)}{\partial a}>\frac 1 2$, proves the uniqueness 
of the pair $(\bar a, t_{int}(a))$. 

As far as the existence is concerned, it is not hard to compute that the two fronts corresponding to positive and negative 
parameters  transversally intersect at the point $(x_{int},y_{int})$. Hence, the jets computed above
are sufficient to show that the geodesic corresponding to $a$ intersects a unique geodesic corresponding to an 
initial condition of the form $\bar a=-a-a^2\psi'(0)+o(a^2)$ at a time of the form $t_{int}(a)=|a|(1 +\frac 1 2 a\psi'(0)+o(a))$.

\hfill $\blacksquare$

\bl\label{conjugate-time}
The conjugate time of a geodesic is strictly bigger than its cut time.
\el
\noindent {\bf Proof:} 
Thanks to the previous computations, 
the absolute value of the Jacobian of the map $(t,a)\mapsto(x(t,a),y(t,a))$ is 
$2 a + 3 a^2 \psi'(0)+ a^3 \Xi(a,\frac t a)$ where $\Xi$ is a smooth function.
This allows to conclude that for $|\frac t a|<2$ and $a$ small enough, the Jacobian is nonzero, 
hence $t$ is not a conjugate time. \hfill $\blacksquare$

\medskip

\noindent{\bf End of the proof.} Let us show that if $a$ is small enough, then the geodesic corresponding to $a$ is optimal on $[0,t_{int}(a)]$. 
By contradiction, assume the geodesic looses optimality at a time  $0<\bar t<t_{int}(a)$. 
\begin{enumerate}
\item $\bar t$ is not a conjugate time thanks to Lemma \ref{conjugate-time}. Hence there exists $b\neq 0$ such that the geodesics corresponding to $a$ and $b$ intersect at time $\bar t$ and they are both optimal on $[0,\bar t]$. 
\item $\bar t$ cannot satisfy $|\frac{\bar t}{b}|\geq 2$ since this would imply that the geodesic corresponding to $b$ is not optimal until $\bar t$ 
thanks to Lemma \ref{calcolo-cut}. 
\item If $|\frac{\bar t}{b}|< 2$, by the uniqueness property given in Lemma~\ref{calcolo-cut}, the parameters $b$ and $a$ cannot have opposite 
signs. Nevertheless they can neither have the same sign as Lemma \ref{derivee} implies that two geodesics corresponding to different parameters of the 
same sign cannot intersect at time $\bar t$.
\end{enumerate}
Hence $t_{int}(a)$ is the cut time of the geodesic corresponding to the parameter $a$ and the cut point is given in \r{cut-sopra}. 
By the transversality of the two fronts (corresponding to initial conditions $a$ and $\bar a$) at time $t_{int}$, the set $\cut_\Zz^+$ 
is locally a 1-dimensional manifold. Moreover, formula \r{cut-sopra} implies that $\cut_\Zz^+$ is  transversal to the distribution at $(0,0)$, 
its tangent vector at $(0,0)$ being $(-\psi'(0)/2, 1)$.

\subsection{The lower part of the cut locus}\label{subs-sotto}

Reasoning as in section~\ref{subs-sopra}, we consider the geodesic starting from $\Zz$, realizing the 
distance from $\Zz$ and entering the lower domain, that is with the initial conditions 
\bqn\label{destra-giu}
\ba{ll}
x(t=0,a)=a , & p_x(t=0,a)=\sign(a),\\
y(t=0,a)= a^2\psi(a), & p_y(t=0,a)=-\displaystyle\frac{\sign(a)}{2a\psi(a)+a^2\psi'(a)}.
\ea
\eqn

\noindent{\bf Computation of the jets of the geodesics.} For $a>0$, setting $\eta=\sqrt a$ and $s=\frac t \eta$, one can check that if $x,y,p_x,p_y$ have orders  
in $\eta$ higher or equal to $1,3,0,-2$, respectively,  then the dynamics has the same or higher orders. As a consequence, since the initial condition respects these orders, one can compute jets with respect to $\eta$ of the solution of system \r{ham-sys-rep} under the form
$$
\begin{array}{lllcclll}
x(s,\eta) & = & \eta x_0(s) + \eta^2 x_1(s) +\eta^3\bar x(s,\eta), &&& 
p_x(s,\eta) & = & {p_x}_0(s) + \eta {p_x}_1(s) +\eta^2\bar p_x(s,\eta),\\
y(s,\eta) & = & \eta^3 y_0(s) + \eta^4 y_1(s) +\eta^5\bar y(s,\eta), &&& 
p_y(s,\eta) & = & \eta^{-2} {p_y}_0(s) +\eta^{-1} {p_y}_1(s) +\bar p_y(s,\eta),
\end{array}
$$
where $\bar x, \bar y, \bar p_x, \bar p_y$ are smooth functions. From the initial conditions \r{destra-giu}, we deduce
$$
\begin{array}{llll}
x_0(0)=0, & x_1(0)=1, & {p_x}_0(0)=1, & {p_x}_1(0)=0,\\
y_0(0)=0, & y_1(0)=1, &
{p_y}_0(0)=-\frac{1}{2}, & {p_y}_1(0)=0,
\end{array}
$$
and using system \r{ham-sys-rep}, the functions $x_0,x_1,y_0,y_1,{p_x}_0,{p_x}_1,{p_y}_0,{p_y}_1$ satisfy
\bqn\label{sys-first-order}
\left\{
\begin{array}{ll}
\dot x_0&={p_x}_0\\
\dot y_0&=\gamma^2{p_y}_0 x_0^4\\
\dot {p_x}_0&=-2\gamma^2{p_y}_0^2x_0^3\\
\dot {p_y}_0&=0
\end{array}\right.\,\,
\left\{
\ba{ll}
 \dot x_1&={p_x}_1\\
\dot y_1&=\gamma^2({p_y}_1x_0^4-2{p_y}_0x_0^2(y_0-2x_0x_1-\alpha x_0^3))\\
\dot {p_x}_1&=\gamma^2{p_y}_0x_0(-4{p_y}_1x_0^2+2{p_y}_0y_0-6{p_y}_0x_0x_1
-5\alpha{p_y}_0x_0^3)\\
 \dot {p_y}_1&=\gamma^2{p_y}_0x_0^2
\ea\right.
\eqn
where $\gamma=\e^{\xi(0,0)}$ and $\alpha=\psi'(0)+\frp{\xi}{x}(0,0)$. 
Thus ${p_y}_0\equiv -\frac{1}{2}$ and one can prove (see   \cite{ag-bonn,BCGJ}) that 
\begin{eqnarray}
x_0(s) & = & -\frac{\sqrt 2}{\sqrt\gamma}\cn(\KKK+\sqrt\gamma s),\nonumber\\
y_0(s) & = & -\frac{2}{3\sqrt\gamma}(\sqrt\gamma s+2\sn(\KKK+\sqrt\gamma s)\cn(\KKK+\sqrt\gamma s)\dn(\KKK+\sqrt\gamma s)),\nonumber
\end{eqnarray}
where $\KKK$ is the complete elliptic integral of the first kind of modulus $\frac{1}{\sqrt 2}$, and $\cn$, $\sn$ and $\dn$ denote the classical Jacobi functions of modulus $\frac{1}{\sqrt 2}$. Recall that
the Jacobi functions $\cn,\sn$ are $4\KKK$-periodic, when $\dn$ is  $2\KKK$-periodic.

Denote by $x_{10},y_{10},{p_x}_{10},{p_y}_{10}$ the solution of the second system in \r{sys-first-order} with $\alpha=0$.  Define $g_1,g_2,g_3,g_4$ by
$$
\ba{llllll}
x_1 & = & x_{10}+\alpha g_1,& {p_x}_1 & = & {p_x}_{10}+\alpha g_3,\\
y_1 & = & y_{10}+\alpha g_2, &{p_y}_1 & = & {p_y}_{10}+\alpha g_4.
\ea
$$
It is easy to see that   the $g_i$ satisfy
\bqn\label{sys-gi}
\left\{
\ba{lll}
g_4 & \equiv &  0\\
\dot g_1 & = & g_3\\
\dot g_2 & = & -\gamma^2x_0^3(2g_1+x_0^2)\\
\dot g_3 & = & -\frac{1}{4}\gamma^2x_0^2(6g_1+5x_0^2),
\ea\right.
\eqn
and the initial conditions are $g_1(0)=g_2(0)=g_3(0)=0$. Notice moreover that, if
 $$
(x_0,y_0,{p_x}_0,{p_y}_0,x_{10},y_{10},{p_x}_{10},{p_y}_{10},g_1,g_2,g_3)
$$ 
is the solution of \r{sys-first-order}, \r{sys-gi}
with initial condition
$(0,0,1,-1/2,1,1,0,0,0,0,0)$ 
then the solution of  \r{sys-first-order}, \r{sys-gi} with initial condition $(0,0,-1,-1/2,-1,1,0,0,0,0,0)$ is
$$
(-x_0,y_0,-{p_x}_0,{p_y}_0,-x_{10},y_{10},-{p_x}_{10},{p_y}_{10},g_1,-g_2,g_3),
$$
which corresponds to the geodesic starting from $\Zz$ with the initial condition $-a<0$.

\medskip


\bl\label{derivee-down}
If $\delta>0$ and $\eta\neq 0$ are small enough and 
$0<|\frac{t}{\eta}|<\frac{2 \KKK}{\sqrt{\gamma}}+\delta$ then $\frp{x(t,a)}{a}>0$. 
\el

%
%

\noindent{\bf Proof:}  Assume $a>0$ (the computation being the same for $a<0$). Then
$$
\frp{x(t,a)}{\eta}=x_0\left(\frac{t}{\eta}\right)-\frac{t}{\eta}\dot x_0\left(\frac{t}{\eta}\right)+\eta\left(2x_1\left(\frac{t}{\eta}\right)
-\frac{t}{\eta}\dot x_1\left(\frac{t}{\eta}\right)\right)+\eta^2x_r\left(\eta,\frac{t}{\eta}\right),
$$
where $x_r$ is a smooth function. Now, the function $f:u\mapsto x_0(u)-u\dot x_0(u)$ is such that $f(0)=0$ and 
$$
f'(u)=-u\ddot x_0(u)=\frac{1}{2}u\gamma^2x_0^3(u)>0 \mbox{ for } u\in\left]0, \frac{2 \KKK}{\sqrt{\gamma}}\right[.
$$ 
Hence, for $\epsilon$ small enough, there exists $\delta>0$ small enough such that $f(u)>\epsilon$ for 
$u\in]\delta,\frac{2 \KKK}{\sqrt{\gamma}}+\delta[$. 
Therefore, if 
$\delta<\frac{t}{\eta}<\frac{2 \KKK}{\sqrt{\gamma}}+\delta$, then $\frp{x(t,a)}{\eta}>\epsilon+\eta\left(2x_1\left(\frac{t}{\eta}\right)
-\frac{t}{\eta}\dot x_1\left(\frac{t}{\eta}\right)\right)+\eta^2x_r\left(\eta,\frac{t}{\eta}\right)>0$ for $\eta$ small enough. 

For $0<\frac{t}{\eta}<\delta$ (possibly reducing  $\delta$ and $\eta$),  since $2x_1(\frac{t}{\eta})-\frac{t}{\eta}\dot x_1(\frac{t}{\eta})=2$ for $t=0$, we have that $\frp{x(t,a)}{\eta}>0$.
\hfill $\blacksquare$

\medskip


\bl\label{calcolo-cut-down}
A geodesic with $a>0$ intersects a geodesic with $\bar a<0$ at time $t_{int}(a)$  where
$$
\ba{l}
\bar a=-a+\frac{\sqrt\gamma \alpha }{\KKK}g_2(\frac{2\KKK}{\sqrt\gamma})a\sqrt a +o(a\sqrt a),\\
t_{int}(a)=\frac{2\KKK}{\sqrt\gamma}
(\sqrt a- \sqrt\gamma\;\; \frac{\alpha g_2(\frac{2\KKK}{\sqrt \gamma})-2x_{10}(\frac{2\KKK}{\sqrt \gamma})}{4\KKK}  a+o(a)).
\ea
$$
\el

\noindent {\bf Proof:} In order to find the expressions given in the statement, we proceed as follows: we fix a time $t_0=\frac{2\KKK}{\sqrt \gamma}\eta_0$ 
and we find two parameters $a>0$ and $\bar a<0$ such that the corresponding geodesics intersect at time $t_0$. Indeed, $t_0$ is a natural candidate 
to approximate the intersection of the geodesics with initial conditions $\eta_0^2$ and $-\eta_0^2$ since $\frac{2\KKK}{\sqrt \gamma}$ 
is the half period of $x_0$.

We look for $a$ and $\bar a$ by setting 
$$
\begin{array}{l}
\eta_+=\sqrt a =\eta_0+c_+\eta_0^2+o(\eta_0^2),\\
\eta_-=\sqrt{|\bar a|}=\eta_0+c_-\eta_0^2+o(\eta_0^2),
\end{array}
$$ 
where $c_+$ and $c_-$ are constants to be found. Remark that this is equivalent to choose 
$t_0=\frac{2\KKK}{\sqrt\gamma}(\eta_+-c_+\eta_+^2+o(\eta_+^2))$ and $\eta_-=\eta_++(c_--c_+)\eta_+^2+o(\eta_+^2)$.
The corresponding geodesic parameterized by $s$ are
\begin{eqnarray}
x_+(s) & = & \eta_+ x_0(s)+\eta_+^2(x_{10}(s)+\alpha g_1(s))+o(\eta_+^2),\nonumber\\
y_+(s) & = & \eta_+^3 y_0(s)+\eta_+^4(y_{10}(s)+\alpha g_2(s))+o(\eta_+^4),\nonumber
\end{eqnarray}
and
\begin{eqnarray}
x_-(s) & = & -\eta_- x_0(s)+\eta_-^2(-x_{10}(s)+\alpha g_1(s))+o(\eta_-^2),\nonumber\\
y_-(s) & = & \eta_-^3 y_0(s)+\eta_-^4(y_{10}(s)-\alpha g_2(s))+o(\eta_-^4).\nonumber
\end{eqnarray}

Let us estimate the geodesic corresponding to $\eta_+$ (resp. $\eta_-$) at $s_+=\frac{t_0}{\eta_+}$ 
(resp. $s_-=\frac{t_0}{\eta_-}$). One computes easily that
\begin{eqnarray}
s_+ & = & \frac{2\KKK}{\sqrt \gamma}(1-c_+\eta_0+o(\eta_0)),\nonumber\\
s_- & = & \frac{2\KKK}{\sqrt \gamma}(1-c_-\eta_0+o(\eta_0)),\nonumber
\end{eqnarray}
and 
\begin{eqnarray}
x_+(s_+) & = & \eta_0^2\left(c_+\frac{2\KKK}{\sqrt \gamma}+x_{10}\left(\frac{2\KKK}{\sqrt \gamma}\right)
+\alpha g_1\left(\frac{2\KKK}{\sqrt \gamma}\right)\right)+o(\eta_0^2),\nonumber\\
y_+(s_+) & = & -\eta_0^3\frac{4\KKK}{3\sqrt\gamma}+\eta_0^4\left(-\frac{4\KKK c_+ }{\sqrt\gamma} 
+y_{10}\left(\frac{2\KKK}{\sqrt \gamma}\right)+\alpha g_2\left(\frac{2\KKK}{\sqrt \gamma}\right)\right)+o(\eta_0^4),\nonumber\\
x_-(s_-) & = & \eta_0^2\left(-c_-\frac{2\KKK}{\sqrt \gamma}-x_{10}\left(\frac{2\KKK}{\sqrt \gamma}\right)
+\alpha g_1\left(\frac{2\KKK}{\sqrt \gamma}\right)\right)+o(\eta_0^2),\nonumber\\
y_-(s_-) & = & -\eta_0^3\frac{4\KKK}{3\sqrt\gamma}+\eta_0^4\left(-\frac{4\KKK c_- }{\sqrt\gamma} 
+y_{10}\left(\frac{2\KKK}{\sqrt \gamma}\right)-\alpha g_2\left(\frac{2\KKK}{\sqrt \gamma}\right)\right)+o(\eta_0^4).\nonumber
\end{eqnarray}
Hence, these two geodesics intersect at time $t_0$ for 
\begin{eqnarray}
c_+ & = & \sqrt\gamma\;\; \frac{\alpha g_2(\frac{2\KKK}{\sqrt \gamma})-2x_{10}(\frac{2\KKK}{\sqrt \gamma})}{4\KKK},\nonumber\\
c_- & = & -\sqrt\gamma\;\; \frac{\alpha g_2(\frac{2\KKK}{\sqrt \gamma})+2x_{10}(\frac{2\KKK}{\sqrt \gamma})}{4\KKK}.\nonumber
\end{eqnarray}
The intersection point is
\begin{equation}\label{intersection}
\ba{rcl}
x_{int}(t_0) & = & \eta_0^2\alpha\;\;\frac{2g_1(\frac{2\KKK}{\sqrt \gamma})+g_2(\frac{2\KKK}{\sqrt \gamma})}{2}
+o(\eta_0^2),\\
y_{int}(t_0) & = & -\eta_0^3 \frac{4\KKK}{3\sqrt \gamma}+o(\eta_0^3),
\ea
\end{equation}
and the intersection time satisfies
\bqn
\frac{t_0}{\eta_+}=\frac{2\KKK}{\sqrt\gamma}(1-c_+\eta_0+o(\eta_0))<\frac{2\KKK}{\sqrt\gamma}+\delta,
\label{inQ}
\eqn
 for $a$ small enough and $\delta$ as in Lemma \ref{derivee-down}. 

One can compute that the two fronts (corresponding to positive and negative initial conditions) are transversal at time
$t_{int}(a)$ at the point $(x_{int},y_{int})$.  Hence the computation with the jets allows to conclude that the two geodesics
intersect.

Notice that numerical computations show that $2g_1(\frac{2\KKK}{\sqrt \gamma})+g_2(\frac{2\KKK}{\sqrt \gamma})\neq 0$.
Moreover, thanks to the assumption in Theorem~\ref{prop-cutz}, $\alpha\neq0$. 
\hfill $\blacksquare$

Lemma \ref{calcolo-cut-down} implies that a geodesic with initial condition $a$ looses optimality at $t\leq t_{int}(a)$
which is less than $\frac{2\KKK}{\sqrt{\gamma}}+\delta$ for $a$ small enough.

\medskip


\bl\label{conj-time-down}
The conjugate time of a geodesic is strictly bigger than its cut time.
\el

\noindent{\bf Proof:} Thanks to the computations on the jets,
one can show that the absolute value of the Jacobian of $(s,\eta)\mapsto(x(s,\eta),y(s,\eta))$ is $\eta^3(J_0(s)+\eta J_1(s)   +\eta^2 J_2(s,\eta))$, 
where $J_0(s)=x_0(s)\dot y_0(s)-3 y_0(s)\dot x_0(s)$, $J_1(0)=-4\sign(\dot x_0(0))$, and $J_2$ is a smooth function. 
It was proven in \cite{BC} that $J_0$  never vanishes between $0$ and $\bar s$ with $\bar s>2\KKK/\sqrt{\gamma}$. 
Moreover $J_1(0)$ has the same sign as the function $J_0$ on the interval $]0,\bar s[$. 
This allows to conclude that for $\delta>0$ small enough, the Jacobian never vanishes on the interval 
$[0,\frac{2\KKK}{\sqrt{\gamma}}+\delta[$ which implies that if $\frac{t}{\eta}<\frac{2\KKK}{\sqrt{\gamma}}+\delta$ 
and $a$ is small enough then $t$ is not a conjugate time for the geodesic with initial condition $a$.\hfill$\blacksquare$

\medskip

\noindent{\bf Variations of $x(\cdot,a)$ until $t_{int}(a)$.}
Let us consider a geodesic with $a>0$.
Since $x_0$ satisfies $\frac{d^2x_0}{ds^2}<0$ on $]0,\frac{2\cal{K}}{\sqrt\gamma}[$, 
$\frac{dx_0}{ds}(0)=1$, $\frac{dx_0}{ds}(\frac{2\cal{K}}{\sqrt\gamma})=-1$, one can prove that $x$ satisfies 
$\frp{x}{s}(s,\eta)\geq 0$ on the interval $[0,\frac{\cal{K}}{\sqrt\gamma}+\epsilon]$ and $\frp{x}{s}(s,\eta)\leq 0$ on the interval 
$[\frac{\cal{K}}{\sqrt\gamma}+\epsilon,\frac{2\cal{K}}{\sqrt\gamma}+\delta]$ for $\eta$ small enough
where $\epsilon$ is a small parameter of order 1 in $\eta$. In particular between $t=0$ and $t=t_{int}(a)$, 
$x$ is first increasing and after decreasing until $t_{int}(a)$.

For $a<0$, one can prove the same way that between $t=0$ and $t=t_{int}(a)$, 
$x$ is first decreasing and after increasing until $t_{int}(a)$.

\medskip

\noindent{\bf Estimation of $x(\cdot,a)$ after the first intersection with the $x$-axis.}
Consider a geodesic with $a>0$. Call $s_a$ the first time $s$ such that $y(a,s)=0$. 
For any $\lambda>0$, one can compute that $y(s,\eta)=\eta^4+o(\eta^4)>0$ 
for $\eta$ small enough and any $s\in[0,\lambda \eta]$.  
Hence $s_a>\lambda\eta$ for $\eta$ small enough.
Since $x(\lambda\eta,\eta)=\eta^2(1+\lambda)+o(\eta^2)$,
fixing $\lambda>\alpha\;\;\frac{2g_1(\frac{2\KKK}{\sqrt \gamma})+g_2(\frac{2\KKK}{\sqrt \gamma})}{2}$,
and thanks to the previous considerations on the variations of $x$, 
the minimum of the $x$-coordinate for $t\in[\eta s_a,t_{int}(a)]$ is attained at $t=t_{int}(a)$.

For $a<0$, one can prove the same way that the maximum of the $x$-coordinate for $t\in[\eta s_a,t_{int}(a)]$ is 
attained at $t=t_{int}(a)$.

Notice that for any geodesic of the 2-ARS the $y$-coordinate is monotone thanks to the normal form (\F3).  Indeed, one easily proves that horizontal lines parameterized by arclength are optimal geodesics. 

\medskip

\noindent{\bf End of the proof.} Now, we have all the ingredients to conclude.

\begin{enumerate}
\item A geodesic cannot loose optimality by reaching its conjugate locus thanks to Lemma \ref{conj-time-down}.

\item Assume that two geodesics with initial conditions $a\neq a'$ of the same sign 
loose optimality by intersecting one each other. They should intersect before $t_{int}(a)$ and $t_{int}(a')$,
i.e., for $t$ such that $0<|\frac{t}{\eta}|<\frac{2 \KKK}{\sqrt{\gamma}}+\delta$ and $0<|\frac{t}{\eta'}|<\frac{2 \KKK}{\sqrt{\gamma}}+\delta$.
Hence Lemma \ref{derivee-down} applies, which leads to a contradiction with $x(t,a)=x(t,a')$.\\
As a consequence a geodesic with parameter $a$ looses optimality by intersecting another geodesic with parameter $a'$ such that
$aa'<0$.

\item Thanks to the monotonicity of the $y$-coordinate, two geodesics with initial conditions $a$ and $a'$ of opposite sign 
can intersect only in the half plane $y<0$. 

In the following, we assume that $\alpha\;\;\frac{2g_1(\frac{2\KKK}{\sqrt \gamma})+g_2(\frac{2\KKK}{\sqrt \gamma})}{2}<0$, i.e., 
$x(t_{int}(a),a)<0$ for $a$ small enough, the proof being the same in the opposite case.

\item Let us consider the geodesics corresponding to $a>0$ and $\bar a<0$ as in Lemma \ref{calcolo-cut-down}. 
For $a'$ such that $0<a'<a$ we have that
$$x(t,a')>x(t_{int}(a'),a')>x(t_{int}(a),a)=x(t_{int}(\bar a),\bar a)>x(t,\bar a)$$
for $t\leq t_{int}(a')$ such that $y(t,\bar a)<0$, where the first and the last inequalities follow from the estimations of $x$ given above. 
This implies that the geodesic corresponding to $a'$ cannot intersect the one corresponding to $\bar a$ before 
loosing optimality.

\item Thanks to Lemma \ref{derivee-down}, if $a'>a$ then $x(t,a)<x(t,a')$ for $t< t_{int}(\bar a)=t_{int}(a)$. Moreover, the estimations of $x$ given above
imply that $x(t,\bar a)<x(t,a)$ for $t< t_{int}(\bar a)$ such that $y(t,a)<0$ and $y(t,\bar a)<0$.
Hence the geodesic corresponding to $a'$
cannot intersect the geodesic corresponding to $\bar a$ before $t_{int}(a)$.

Finally we conclude that any geodesic with parameter $\bar a<0$ looses optimality by intersecting the geodesic with parameter
$a>0$, which implies that the geodesic with parameter $a$ looses optimality at the same time.
\end{enumerate}

Since the two fronts are transversal at the intersection of the two geodesics corresponding to the initial conditions $a$ and $\bar a$, 
the lower part of the cut locus is locally a 1-dimensional manifold. Together with the formulae (\ref{intersection}), this implies that the 
lower part of the cut is half a cusp tangent to $y=0$. The set $\cut_\Zz^-$ for an example of 2-ARS is portraited in Figure~\ref{f-lower}.

 \begin{figure}[h!]
\begin{center}
\input{front2.pstex_t}
\caption{The set $\cut_\Zz^-$ for the ARS with orthonormal frame $F_1=\frp{}{x},\, F_2=(y-x^2-x^3)\frp{}
{y}$ }\label{f-lower}
\end{center}
\end{figure}

 \brem
In \cite{BCGJ}  the description of the cut locus to a tangency point was given. In that paper, only the existence of the intersection point  is actually proved.
The same arguments as in Section~\ref{subs-sotto}
can be  used to prove that the intersection really corresponds to the cut locus.
\erem

\section{Proof of Theorem~\ref{rnf-tangency} }\label{s-6}

Let us start by proving the existence of a subset of $\spadesuit$ passing through the tangency point $q$ and satisfying the following conditions: {\bf i)} 
it is the support of a smooth curve,  {\bf ii)} it has a tangent direction which is transversal to the distribution at $q$.

Choose a local representation of the type (\F3). 
By construction, $K$ is well defined outside the singular set $\Zz$. 
 The set $\spadesuit\setminus\Zz$ is implicitly defined by the equation
\begin{equation}\label{eqimpl}
G(\nabla(||\nabla K||^2),(\nabla K)^\perp)=0.
\end{equation}
Computing the left hand side of equation \r{eqimpl}, using the expression of the curvature (\ref{curvature}), we find that
$$
G(\nabla(||\nabla K||^2),(\nabla K)^\perp)=\frac{e^{2\xi(x,y)}h(x,y)}{(y-x^2\psi(x))^{8}},
$$
where $h$ is a smooth function. 
Hence, equation \r{eqimpl} is equivalent to $h(x,y)=0$. The development of $h$ at the point $(0,0)$ is 
$$
h(x,y)=\omega\big( y^4(10 \psi(0)^2 x + y(3\psi'(0)-2\partial_x\xi(0,0)\psi(0)))+ \sum_{i=0}^6a_i(x,y)x^iy^{6-i}\big),
$$
where $\omega$ is a nonzero constant and  $a_i$ are smooth functions.
Let us show that there exists a smooth function $b:I\rightarrow \R$ defined on a neighborhood $I$ of 0 such that after the coordinate change 
\begin{equation*}
\bar x=10 \psi(0)^2 x + y(3\psi'(0)-2\partial_x\xi(0,0)\psi(0))-b(y)y^2,\quad \bar y=y,
\end{equation*}
 we have
$h(x(\bar x,\bar y),y(\bar x,\bar y))=\bar x\overline h(\bar x,\bar y)$. 
In the new coordinate system, we have
$$
h(x(\bar x,\bar y),y(\bar x,\bar y))=\omega(\bar y^4\bar x+ F(\bar x,\bar y)), \mbox{ where }
F(\bar x,\bar y)=b(\bar y)\bar y^6+
 \sum_{i=0}^6a_i(x(\bar x,\bar y),\bar y)
(x(\bar x,\bar y))^i\bar y^{6-i}.
$$
In order $\bar x$ to be factorizable in $F$, we require that $F(0,\bar y)\equiv 0$. Since $F(0,\bar y)=\bar y^6 R(b(\bar y),\bar y)$, where
$$
R(b(\bar y),\bar y)=b(\bar y)+ \sum_{i=0}^6 \frac{a_i(x(0,\bar y),\bar y)}{10^i\psi(0)^{2 i}}
(-3\psi'(0)+2\partial_x\xi(0,0)\psi(0)-b(\bar y)\bar y)^i,
$$ 
it follows that $F(0,\bar y)\equiv 0$ if and only if there exists  a smooth function $b$ defined on a neighborhood of 
$0$ such that $R(b(\bar y),\bar y)\equiv 0$.
Let $\overline b=- \sum_{i=0}^6 \frac{a_i(0,0)}{10^i\psi(0)^{2 i}}
(-3\psi'(0)+2\partial_x\xi(0,0)\psi(0))^i$. Then, 
since $(b,\bar y)\mapsto R(b,\bar y)$ is smooth, $R(\overline b,0)=0$, $\partial_bR(\overline b,0)=1$,  by the implicit function theorem there exists a smooth function $b(\bar y)$ with the properties above. Therefore, coming back to the $(x,y)$ coordinates 
we have shown that 
$$
h(x,y)=\omega (10 x + y (3\psi'(0)-2\partial_x\xi(0,0))+b(y)y^2) (y^4+\tilde F(x,y)),
$$
where $\tilde F$ is smooth function of order 5 in $(x,y)$ and $b$ is the function built  above.
The last equation implies that the set $C=\{(x,y)\mid 10 x + y (3\psi'(0)-2\partial_x\xi(0,0))+b(y)y^2=0 \}$ is a connected component of the set $\spadesuit$, it passes through $(0,0)$, it is smooth at $(0,0)$ and its tangent line at $(0,0)$ is 
$$x=\frac{1}{10 \psi(0)^2}(2 \partial_x\xi(0,0)\psi(0)-3 \psi'(0))y,$$ that is transversal to the distribution at $(0,0)$. 

Moreover, since $\tilde F$ has order 5 in $(x,y)$,  any curve contained in the set $\{(x,y)\mid h(x,y)=0\}$ but not
in $C$ must have a tangent line at $(0,0)$ belonging to the distribution at $(0,0)$.  

Requiring that the tangent direction to $C$ at $(0,0)$ is vertical we get the condition 
\begin{equation}\label{eq-condtang}
2 \partial_x\xi(0,0)\psi(0)-3 \psi'(0)=0.
\end{equation}
As a consequence, any curve transversal to $\Zz$ at $q$ contained in $\spadesuit$ should have $C$ as support.

\medskip

Now, let us prove that there is a canonical parametrization of $C$. Let us choose any parameterization of $C$
and construct $(x,y)$ and $(X,Y)$ as in Procedure 1. 
Then, since $\dim(\bD(q))=\dim\bD_2(q)=1$ and $\dim\bD_3(q)=2$, one gets that $Y(0,0)$ and $[X,Y](0,0)$ are horizontal and $[X,[X,Y]](0,0)$ is not horizontal, 
which is equivalent to $f(0,0)=f_x(0,0)=0$ and $f_{xx}(0,0)\neq 0$.

Assume $f_{xx}(0,0)>0$. Denote by $c(\cdot)$ the parameterized curve whose support is $C$, such that $c(0)=q$, and satisfying $\dot c(y)=1/2f_{xx}(0,y)\partial_y$ in the $(x,y)$ coordinates.  With this choice for the parameterization, if $(\bar x, \bar y)$ and $(\bar X, \bar Y)$ denotes the coordinates and the local representation defined by $c(\cdot)$ via Procedure 1, then $[\bar X,[\bar X,\bar Y]](0,\bar y)=2 \partial_{\bar y}$. Equivalently, if $\bar Y=\bar f(\bar x, \bar y)\partial_{\bar y}$ we have $\bar f_{\bar x \bar x}(0,\bar y)\equiv 2$. An easy computation shows that choosing the parameterization $y\mapsto c(-y)$, the condition $\bar f_{\bar x \bar x}(0,\bar y)\equiv 2$ is still fulfilled.

If $f_{xx}(0,0)<0$, the same arguments prove that, denoting by $c(\cdot)$ the curve whose support is $C$, such that $c(0)=q$ and satisfying $\dot c(y)=-1/2 f_{xx}(0,y)\partial_y$ in the $(x,y)$ coordinates, then $\bar f_{\bar x \bar x}(0,\bar y)\equiv-2$. Moreover, the last condition does not change when choosing the opposite orientation for $c(\cdot)$.

The curve $c(\cdot)$ defined above satisfies (i), (ii) and, depending on the sign of $\bar f_{\bar x\bar x}(0,0)$, the first part of (iii). It is uniquely determined
up to orientation. One can prove easily that if 
$\bar f_{\bar y}\partial_{\bar y}=[\partial_{\bar y},\bar Y]=\lambda\partial_{\bar y}$ then, changing the parameterization
of $c(\cdot)$ for $y\mapsto c(-y)$ one gets $\bar f_{\bar y}\partial_{\bar y}=[\partial_{\bar y},\bar Y]=-\lambda\partial_{\bar y}$.

In Theorem \ref{rnf-tangency}, we fix the orientation of $c(\cdot)$ (and hence we fix $c(\cdot)$) by asking that when 
$f_{\bar x\bar x}(0,0)>0$ then $\bar f_{\bar y}(0,0)<0$ and when
$f_{\bar x\bar x}(0,0)<0$ then $\bar f_{\bar y}(0,0)>0$. It is not hard to prove that the first situation corresponds to a point of type \tm\ and the
second to a point of type \tp.

\section{Proof of Corollary~\ref{main-theorem}}\label{s-7}

\noindent{\bf Grushin points.}
Equation (Ga) is a consequence of the fact that the vertical axis $x=0$ is contained in
$\Zz$ and hence the distribution has dimension one for $x=0$. 
Equation (Gb) expresses the fact that $[X_1,X_2]=\partial_y$ along 
the vertical axis, by construction.

\medskip

\noindent{\bf Tangency points.}
Denote by $(x,y)$ the coordinate system 
and by $(X,Y)$ the orthonormal frame given by Procedure 1 when the curve $c(\cdot)$ is the one given by Theorem~\ref{rnf-tangency}. 
The fact that $\Delta(q)$ and $\Delta_2(q)$ have dimension one implies (Ta) and (Tb). 
Equations (Tc) and (Td) are implied by (iii) of Theorem~\ref{rnf-tangency}. Equation (Te) is a direct consequence of (ii)  of Theorem~\ref{rnf-tangency}.

\medskip

\noindent{\bf Riemannian points of type 1.}
Properties (R1a), (R1b) and (R1c) are direct consequences of the discussion given in Section \ref{tcR1}.

Equation (R1d) comes from property (R1a). 
For what concerns (R1e) and (R1f), let us compute $\nabla K$ using the equality (\ref{curvature})
 $$
 \nabla K(x,y)=\big(\partial_x^3\phi (x,y)-2 \partial_x\phi (x,y) \partial_x^2\phi
   (x,y),e^{2 \phi (x,y)} 
   \left(-2 \partial_x\phi (x,y) \partial_{x}\partial_y\phi
   (x,y)+\partial_x^2\partial_y\phi
   (x,y)\right)\big).
 $$
Applying (R1b) and (R1c) to this equality, one gets (R1e) and (R1f).

\medskip

\noindent{\bf Riemannian points of type 2.}
The properties (R2a) and (R2b) are direct consequences of the discussion given in Section \ref{tcR2}.
Inequality (R2c) follows from the fact that the vertical axis is a crest and the horizontal one an anticrest.
Inequality (R2d)  is implied by the fact that the vertical axis is a valley and the horizontal one an antivalley.
Inequality (R2e) is a consequence of the fact  that the vertical axis is a crest and the horizontal one a valley.
Inequality (R2f) fixes the orientation of the parameterization (cf. Section \ref{tcR2}).

\bibliographystyle{abbrv}
\bibliography{biblio_AR-invariants}

\begin{thebibliography}{10}

\bibitem{agrompatto}
A.~Agrachev.
\newblock Compactness for sub-{R}iemannian length-minimizers and
  subanalyticity.
\newblock {\em Rend. Sem. Mat. Univ. Politec. Torino}, 56(4):1--12 (2001),
  1998.
\newblock Control theory and its applications (Grado, 1998).

\bibitem{AgrBarBoscbook}
A.~Agrachev, D.~Barilari, and U.~Boscain.
\newblock {\em Introduction to Riemannian and sub-Riemannian geometry (Lecture
  Notes)}.
\newblock http://people.sissa.it/agrachev/agrachev\_files/notes.html.

\bibitem{ag-bonn}
A.~Agrachev, B.~Bonnard, M.~Chyba, and I.~Kupka.
\newblock Sub-{R}iemannian sphere in {M}artinet flat case.
\newblock {\em ESAIM Control Optim. Calc. Var.}, 2:377--448, 1997.

\bibitem{ABS}
A.~Agrachev, U.~Boscain, and M.~Sigalotti.
\newblock A {G}auss-{B}onnet-like formula on two-dimensional
  almost-{R}iemannian manifolds.
\newblock {\em Discrete Contin. Dyn. Syst.}, 20(4):801--822, 2008.

\bibitem{zelclass}
A.~Agrachev and I.~Zelenko.
\newblock On feedback classification of control-affine systems with one- and
  two-dimensional inputs.
\newblock {\em SIAM J. Control Optim.}, 46(4):1431--1460 (electronic), 2007.

\bibitem{agra-gauss}
A.~A. Agrach{\"e}v.
\newblock A ``{G}auss-{B}onnet formula'' for contact sub-{R}iemannian
  manifolds.
\newblock {\em Dokl. Akad. Nauk}, 381(5):583--585, 2001.

\bibitem{euler}
A.~A. Agrachev, U.~Boscain, G.~Charlot, R.~Ghezzi, and M.~Sigalotti.
\newblock Two-dimensional almost-{R}iemannian structures with tangency points.
\newblock {\em Ann. Inst. H. Poincar\'e Anal. Non Lin\'eaire}, 27(3):793--807,
  2010.

\bibitem{book2}
A.~A. Agrachev and Y.~L. Sachkov.
\newblock {\em Control theory from the geometric viewpoint}, volume~87 of {\em
  Encyclopaedia of Mathematical Sciences}.
\newblock Springer-Verlag, Berlin, 2004.
\newblock Control Theory and Optimization, II.

\bibitem{baouendi}
M.~S. Baouendi.
\newblock Sur une classe d'op\'erateurs elliptiques d\'eg\'en\'er\'es.
\newblock {\em Bull. Soc. Math. France}, 95:45--87, 1967.

\bibitem{bellaiche}
A.~Bella{\"{\i}}che.
\newblock The tangent space in sub-{R}iemannian geometry.
\newblock In {\em Sub-{R}iemannian geometry}, volume 144 of {\em Progr. Math.},
  pages 1--78. Birkh\"auser, Basel, 1996.

\bibitem{BCa}
B.~Bonnard and J.~B. Caillau.
\newblock {S}ingular {M}etrics on the {T}wo-{S}phere in {S}pace {M}echanics.
\newblock Preprint 2008, HAL, vol. 00319299, pp. 1-25.

\bibitem{tannaka}
B.~Bonnard, J.-B. Caillau, R.~Sinclair, and M.~Tanaka.
\newblock Conjugate and cut loci of a two-sphere of revolution with application
  to optimal control.
\newblock {\em Ann. Inst. H. Poincar\'e Anal. Non Lin\'eaire},
  26(4):1081--1098, 2009.

\bibitem{BCGJ}
B.~Bonnard, G.~Charlot, R.~Ghezzi, and G.~Janin.
\newblock {T}he {S}phere and the {C}ut {L}ocus at a {T}angency {P}oint in
  {T}wo-{D}imensional {A}lmost-{R}iemannian {G}eometry.
\newblock {\em J. Dynam. Control Systems}, 17(1):141--161, 2011.

\bibitem{BC}
B.~Bonnard and M.~Chyba.
\newblock M\'ethodes g\'eom\'etriques et analytiques pour \'etudier
  l'application exponentielle, la sph\`ere et le front d'onde en g\'eom\'etrie
  sous-riemannienne dans le cas {M}artinet.
\newblock {\em ESAIM Control Optim. Calc. Var.}, 4:245--334 (electronic), 1999.

\bibitem{q4}
U.~Boscain, T.~Chambrion, and G.~Charlot.
\newblock Nonisotropic 3-level quantum systems: complete solutions for minimum
  time and minimum energy.
\newblock {\em Discrete Contin. Dyn. Syst. Ser. B}, 5(4):957--990, 2005.

\bibitem{BCha}
U.~Boscain and G.~Charlot.
\newblock Resonance of minimizers for {$n$}-level quantum systems with an
  arbitrary cost.
\newblock {\em ESAIM Control Optim. Calc. Var.}, 10(4):593--614 (electronic),
  2004.

\bibitem{q1}
U.~Boscain, G.~Charlot, J.-P. Gauthier, S.~Gu{\'e}rin, and H.-R. Jauslin.
\newblock Optimal control in laser-induced population transfer for two- and
  three-level quantum systems.
\newblock {\em J. Math. Phys.}, 43(5):2107--2132, 2002.

\bibitem{BCGS}
U.~Boscain, G.~Charlot, R.~Ghezzi, and M.~Sigalotti.
\newblock Lipschitz classification of almost-riemannian distances on compact
  oriented surfaces.
\newblock {\em Journal of Geometric Analysis}, pages 1--18.
\newblock 10.1007/s12220-011-9262-4.

\bibitem{camillo}
U.~Boscain and C.~Laurent.
\newblock {T}he {L}aplace--{B}eltrami operator in almost-{R}iemannian
  {G}eometry.
\newblock Preprint 2011, arXiv:1105.4687.

\bibitem{high-order}
U.~Boscain and M.~Sigalotti.
\newblock High-order angles in almost-{R}iemannian geometry.
\newblock In {\em Actes de {S}\'eminaire de {T}h\'eorie {S}pectrale et
  {G}\'eom\'etrie. {V}ol. 24. {A}nn\'ee 2005--2006}, volume~25 of {\em S\'emin.
  Th\'eor. Spectr. G\'eom.}, pages 41--54. Univ. Grenoble I, 2008.

\bibitem{FL1}
B.~Franchi and E.~Lanconelli.
\newblock Une m\'etrique associ\'ee \`a une classe d'op\'erateurs elliptiques
  d\'eg\'en\'er\'es.
\newblock {\em Rend. Sem. Mat. Univ. Politec. Torino}, (Special Issue):105--114
  (1984), 1983.
\newblock Conference on linear partial and pseudodifferential operators
  (Torino, 1982).

\bibitem{grushin1}
V.~V. Gru{\v{s}}in.
\newblock A certain class of hypoelliptic operators.
\newblock {\em Mat. Sb. (N.S.)}, 83 (125):456--473, 1970.

\bibitem{jean1}
F.~Jean.
\newblock Uniform estimation of sub-{R}iemannian balls.
\newblock {\em J. Dynam. Control Systems}, 7(4):473--500, 2001.

\bibitem{kulkarni}
R.~S. Kulkarni.
\newblock Curvature and metric.
\newblock {\em Ann. of Math. (2)}, 91:311--331, 1970.

\bibitem{pelletier}
F.~Pelletier.
\newblock Quelques propri\'et\'es g\'eom\'etriques des vari\'et\'es
  pseudo-riemanniennes singuli\`eres.
\newblock {\em Ann. Fac. Sci. Toulouse Math. (6)}, 4(1):87--199, 1995.

\bibitem{pelle2}
F.~Pelletier and L.~Val{\`e}re~Bouche.
\newblock The problem of geodesics, intrinsic derivation and the use of control
  theory in singular sub-{R}iemannian geometry.
\newblock In {\em Actes de la {T}able {R}onde de {G}\'eom\'etrie
  {D}iff\'erentielle ({L}uminy, 1992)}, volume~1 of {\em S\'emin. Congr.},
  pages 453--512. Soc. Math. France, Paris, 1996.

\bibitem{pontryagin-book}
L.~S. Pontryagin, V.~G. Boltyanski{\u\i}, R.~V. Gamkrelidze, and E.~F.
  Mishchenko.
\newblock {\em The {M}athematical {T}heory of {O}ptimal {P}rocesses}.
\newblock ``Nauka'', Moscow, fourth edition, 1983.

\bibitem{jean2}
M.~Vendittelli, G.~Oriolo, F.~Jean, and J.-P. Laumond.
\newblock Nonhomogeneous nilpotent approximations for nonholonomic systems with
  singularities.
\newblock {\em IEEE Trans. Automat. Control}, 49(2):261--266, 2004.

\end{thebibliography}

\end{document}